\RequirePackage{fix-cm}

\documentclass[smallcondensed]{svjour3}     
\smartqed  
\usepackage{graphicx}
\usepackage{amsmath,bm,amssymb}
\usepackage{color}
\newcommand{\eps}{\varepsilon}
\numberwithin{equation}{section}
\numberwithin{theorem}{section}
\numberwithin{lemma}{section}
\numberwithin{remark}{section}
\graphicspath{{figures/}}
\newcommand{\nn}{\nonumber}

\begin{document}

\title{Improved uniform error bounds of the time-splitting methods for the long-time  (nonlinear) Schr\"odinger equation}

\titlerunning{Improved uniform error bounds for time-splitting methods for NLSE}   

\author{Weizhu Bao \and
Yongyong Cai \and
Yue Feng }

\authorrunning{W. Bao, Y. Cai and Y. Feng}

\institute{W. Bao \at
              Department of Mathematics,
National University of Singapore, Singapore 119076 \\
              \email{matbaowz@nus.edu.sg}
              \and
Y. Cai \at
Laboratory of Mathematics and Complex Systems and School of Mathematical Sciences, Beijing Normal University, Beijing 100875, China \\
\email{yongyong.cai@bnu.edu.cn}
           \and
           Y. Feng \at
           Department of Mathematics,
National University of Singapore, Singapore 119076\\
\email{fengyue@u.nus.edu}
}

\date{Received: date / Accepted: date}

\maketitle



\begin{abstract}
We establish improved uniform error bounds for the time-splitting methods for the long-time dynamics of the Schr\"odinger equation with small potential and the nonlinear Schr\"odinger equation (NLSE) with weak nonlinearity. For the Schr\"odinger equation with small potential characterized by a dimensionless parameter $\varepsilon \in (0, 1]$, we employ the unitary flow property of  the (second-order) time-splitting Fourier pseudospectral (TSFP) method in $L^2$-norm to prove a uniform error bound  at time $t_\varepsilon=t/\varepsilon$ as $C(t)\widetilde{C}(T)(h^m +\tau^2)$ up to $t_\varepsilon\leq T_\varepsilon = T/\varepsilon$ for any $T>0$ and uniformly for $\varepsilon\in(0,1]$, while $h$ is the mesh size, $\tau$ is the time step, $m \ge 2$ and $\tilde{C}(T)$ (the local error bound) depend on the regularity of the exact solution, and $C(t) =C_0+C_1t$ grows at most linearly with respect to $t$ with $C_0$ and $C_1$ two positive constants independent of $T$, $\varepsilon$, $h$ and $\tau$. Then by introducing a new technique of {\sl regularity compensation oscillation} (RCO) in which the high frequency modes are controlled by regularity and the low frequency modes are analyzed by phase cancellation and energy method, an improved uniform (w.r.t $\varepsilon$) error bound at $O(h^{m-1} + \varepsilon \tau^2)$ is established in $H^1$-norm for the long-time dynamics up to the time at $O(1/\varepsilon)$ of the Schr\"odinger equation with $O(\varepsilon)$-potential with $m \geq 3$. Moreover, the RCO technique is extended to prove an improved uniform error bound at $O(h^{m-1} + \varepsilon^2\tau^2)$ in $H^1$-norm for the long-time dynamics up to the time at $O(1/\varepsilon^2)$ of the cubic NLSE with $O(\varepsilon^2)$-nonlinearity strength.   Extensions to
the first-order and fourth-order time-splitting methods are discussed.
 Numerical results are reported to validate our error estimates and to demonstrate that they are sharp.
\end{abstract}

\section{Introduction}
The (nonlinear) Schr\"odinger equation arises in various physical phenomena, such as quantum mechanics, Bose-Einstein condensates, laser beam propagation, plasma and particle physics  \cite{ABB,Caz,ESY,Faou,SS}. In this paper, we consider the following Schr\"odinger equation
\begin{equation}
i\partial_t \psi({\bf x}, t) = - \Delta \psi({\bf x}, t) + \varepsilon V({\bf x}) \psi({\bf x}, t), \ {\bf x} \in \Omega, \ t > 0,
\label{eq:linear}
\end{equation}
and the nonlinear Schr\"odinger equation (NLSE)
\begin{equation}
i\partial_t \psi({\bf x}, t) = - \Delta \psi({\bf x}, t) \pm  \eps^2 \vert\psi({\bf x}, t)\vert^2\psi({\bf x}, t), \ {\bf x} \in \Omega, \ t > 0,
\label{eq:WNE}
\end{equation}
with the initial data
\begin{equation}
\psi({\bf x}, 0) = 	 \psi_0 ({\bf x}), \quad {\bf x} \in  \overline{\Omega},
\label{eq:initial}
\end{equation}
where $\Omega =  \prod_{i = 1}^{d} (a_i, b_i) \subset \mathbb{R}^d$ $(d = 1, 2, 3)$ is a bounded domain equipped with periodic boundary conditions.
Here, $t$ is time, ${\bf x} = (x_1, \cdots, x_d)^T \in \mathbb{R}^d$ is the spatial coordinate, $\psi({\bf x}, t) \in \mathbb{C}$ is the complex order parameter/wave function, $V({\bf x}) \in \mathbb{R}$ is a given external potential, $\varepsilon \in (0, 1]$ is a dimensionless parameter. 
In the Schr\"odinger equation \eqref{eq:linear}, the amplitude of the potential is characterized by the  parameter $\varepsilon \in (0, 1]$. In the NLSE \eqref{eq:WNE}, the strength of the nonlinearity is $O(\eps^2)$  --
{\sl NLSE with weak nonlinearity} --  and the dynamics of the NLSE \eqref{eq:WNE} with $O(1)$-initial data is equivalent to the NLSE with $O(1)$-nonlinearity and $O(\eps)$-initial data -- {\sl NLSE with small initial data}, e.g. by setting $\phi({\bf x},t)=\eps\psi({\bf x},t)$, the  NLSE \eqref{eq:WNE}  with \eqref{eq:initial} becomes
\begin{equation}
\begin{cases}
i\partial_t \phi({\bf x}, t) = - \Delta \phi({\bf x}, t) \pm  \vert\phi({\bf x}, t)\vert^2\phi({\bf x}, t), \ {\bf x} \in \Omega, \ t > 0,\\
\phi({\bf x}, 0) =\varepsilon \psi_0 ({\bf x}):=\phi_0({\bf x})=O(\varepsilon), \quad {\bf x} \in  \overline{\Omega}.
\end{cases}
\end{equation}

 In the past two decades, many accurate and efficient numerical methods have been proposed and analyzed to simulate the (nonlinear) Schr\"odinger equation including the finite difference time domain (FDTD) methods \cite{AK,BC,DFP}, the exponential wave integrator Fourier pseudospectral (EWI-FP) method \cite{CCO,Deb,HO}, the time-splitting Fourier pseudospectral (TSFP) method \cite{BBD,BJP,LU,WH}, etc. Among these numerical methods, the TSFP method  preserves a set of geometric properties and performs much better than the other numerical approaches regarding the stability, efficiency, accuracy and spatial/temporal resolution \cite{ABB,BJP,BJP2}.  However, convergence analysis for the TSFP method applied to the (nonlinear) Schr\"odinger equation is normally valid up to the finite-time dynamics at $O(1)$ and we refer to \cite{BJP,BBD,Faou,HLW,LU,Thal2} and references therein.

Recently, long-time behaviors of the (nonlinear) Schr\"odinger equation on compact domains have received a great deal of attention \cite{BGHS,CCMM,Faou,FGH,FGL}. Along the analytical front, the existence of the solution, the asymptotic behavior and conservation laws have been well studied in the literature \cite{Bour,BGT,CK,HTT,Tao}. From the viewpoint of numerical analysis, the stability of the plane wave solutions and long-time preservations of the actions and energy for the TSFP method have been shown for the NLSE with the help of Birkhoff normal form and the modulated Fourier expansion \cite{CHL,FGL,FGP,GL1,GL2}. For the long-time error estimates of the numerical schemes, improved error bounds for time-splitting methods have been proven under the constraint that the time step $\tau$ is an integer fraction of the period of the principal linear part \cite{CMTZ}. Due to the usage of the properties for the periodic function, extensions of the improved error bounds to higher dimensions require that the aspect ratio of the domain is rational. In addition, error estimates of the splitting methods have been established  with the error bound growing linearly in time  for the Maxwell's equations  \cite{CLL,CLL2} and the Schr\"odinger equations \cite{JL00} (semi-discrete-in-time case). However, such  linear growth of  the fully discrete TSFP error bound  for  the Schr\"odinger equation has not been reported.

The aim of this work is to establish the improved uniform error bounds for the TSFP method for the long-time dynamics of the Schr\"odinger equation with small potential and the NLSE with weak nonlinearity, removing the previous assumptions on the periodicity of the free Schr\"odinger evolutionary operator and the integer fraction time steps. First, we prove a uniform error  bound in $L^2$-norm for the TSFP method applied to the Sch\"odinger equation with the constant in the error bound growing linearly with respect to the time $t$. Based on this error bound, for a given accuracy tolerance $\delta_0$ and time step $\tau$, we could obtain the computational time within the accuracy $\delta_0$ by using the TSFP method is $O(\delta_0/\tau^2)$ for $\eps = 1$, i.e., with the smaller time step $\tau$, the longer dynamics for the Schr\"odinger equation can be calculated! Then by introducing a new technique of {\bf regularity compensation oscillation} (RCO) in which the high frequency modes are controlled by regularity and the low frequency modes are analyzed by phase cancellation and energy method, an improved uniform error bound in $H^1$-norm for the Schr\"odinger equation with $O(\varepsilon)$-potential up to the time $O(1/\varepsilon)$ is carried out at $O(h^{m-1} + \varepsilon\tau^2 + \tau_0^{m-1})$  with $m \geq 3$ depending on the regularity of the exact solution and $\tau_0 \in (0, 1)$ a parameter fixed. In addition, the technique of RCO is extended to the proof of an improved uniform error bound for the cubic NLSE with $O(\varepsilon^2)$-nonlinearity strength up to the time at $O(1/\eps^2)$ with the error bound in $H^1$-norm at $O(h^{m-1} + \varepsilon^2 \tau^2 + \tau_0^{m-1})$.

Here, we briefly explain the idea of our analysis. For sufficiently regular solution, we use the smoothness of the exact solution to control the high frequency modes ($>1/\tau_0$) as $\tau_0^{m-1}$, where $\tau_0$ is a chosen frequency cut-off parameter. The low frequency modes ($\leq1/\tau_0$) will be treated by the RCO technique for sufficiently small $\tau$ and non-resonant $\tau$, which basically asserts that the error of the low frequency part behaves much better (satisfies the improved error bounds) as long as the time step size $\tau$ is non-resonant or resolves the frequency. The regularity  compensation oscillation (RCO) comes from the facts   that the high modes are bounded by the regularity of the exact solution and
an order of $\varepsilon$ could be gained by noticing that $i\partial_t\psi+\Delta \psi=O(\varepsilon\psi)$ from \eqref{eq:linear}, and  respectively, an order of $\varepsilon^2$ could be gained by $i\partial_t\psi+\Delta \psi=O(\varepsilon^2\vert\psi\vert^2\psi)$ from \eqref{eq:WNE}, i.e., a suitable combination of higher order derivatives could compensate the wave oscillation of magnitude $O(\varepsilon)$/$O(\varepsilon^2)$.

The rest of this paper is organized as follows. In Section 2, the uniform error bound for the TSFP method in $L^2$-norm for the Schr\"odinger equation with $O(\varepsilon)$-potential  up to the final time $T_\eps = T/\eps$ is proven and the error is shown to grow linearly  with respect to  $T$.  Then, the improved uniform error bound in $H^1$-norm is rigorously established with the help of a new technique of regularity compensation oscillation (RCO). In Section 3, the RCO technique is extended to  analyze the  improved uniform error bound in $H^1$-norm for the cubic NLSE with $O(\varepsilon^2)$-nonlinearity strength up to the final long-time at $O(1/\varepsilon^2)$. In Sections 2 \& 3, extensive numerical results are reported to validate our error estimates and demonstrate that they are sharp. Finally, some conclusions are drawn in Section 4. Throughout the paper,  the notation $A \lesssim B$ is used to represent that there exists a generic constant $C > 0$, which is independent of the mesh size $h$, time step $\tau$ and $\varepsilon$ such that $\vert A \vert\leq C B$.

\section{Improved uniform error bounds for the Schr\"odinger equation}
In this section, we adopt the time-splitting Fourier pseudospectral (TSFP) method to numerically solve the Schr\"odinger equation \eqref{eq:linear} and rigorously establish the uniform error bound in $L^2$-norm and improved uniform error bound in $H^1$-norm using RCO. For the simplicity of presentation, we only carry out the analysis in one dimension (1D) and generalizations to higher dimensions are straightforward (see also Remark \ref{remark:hd} for discussion). In 1D, the  Schr\"odinger equation \eqref{eq:linear} with the initial data \eqref{eq:initial} and periodic boundary conditions on the domain $\Omega = (a, b)$ can be written as
\begin{equation}
\left\{
\begin{split}
& i\partial_t \psi(x, t) = -\Delta \psi(x, t) + \varepsilon V(x) \psi(x, t),\quad a<x<b,\ t > 0, \\
&\psi(a, t) = \psi(b, t),\ \partial_x \psi(a, t) = \partial_x \psi(b, t),\quad t \geq 0,\\
&\psi(x, 0) = \psi_0(x), \quad x \in [a, b].
\end{split}\right.
\label{eq:linear_1D}
\end{equation}

\subsection{The TSFP method}
By the splitting technique \cite{LU,MQ,Thal}, the  Schr\"odinger equation \eqref{eq:linear_1D} can be decomposed into two subproblems. The first one is
\begin{equation}
\left\{
\begin{split}
& i\partial_t \psi(x, t) = -\Delta  \psi(x, t),\quad x \in \Omega, \quad t > 0, \\
&\psi(a, t) = \psi(b, t),\ \partial_x \psi(a, t) = \partial_x \psi(b, t),\ t \geq 0, \\
& \psi(x, 0) = \psi_0(x), \ x \in [a, b],
\end{split}\right.
\label{eq:sub1}
\end{equation}
which can be solved exactly in phase space
\begin{equation}
\psi(\cdot, t) = e^{it\Delta} \psi_0(\cdot),\quad t \geq 0.
\end{equation}
The second one is to solve
\begin{equation}
\left\{
\begin{split}
& i\partial_t \psi(x, t) = \eps V(x)\psi(x, t),\quad x \in \Omega, \quad t > 0, \\
& \psi(x, 0) = \psi_0(x), \quad x \in [a, b],
\end{split}\right.
\label{eq:sub2}
\end{equation}
which can be integrated exactly in time, for $x \in [a, b]$, as
\begin{equation}
\psi(x, t) = e^{-i \eps t V(x)}\psi_0(x),\quad   t \geq 0.
\end{equation}

Choose $\tau>0$ as the time step size and $t_n=n\tau$ for $n=0, 1,\ldots$ as the time steps. Denote $\psi^{[n]}(x)$ to be the approximation of $\psi(x, t_n)$ for $n \ge 0$, then a second-order semi-discretization of the Schr\"odinger equation \eqref{eq:linear_1D} via the Strang splitting can be given as:
\begin{equation}
\label{eq:Strang}
\psi^{[n+1]}(x) = \mathcal{S}_{\tau}(\psi^{[n]})= e^{i\frac{\tau}{2}\Delta}e^{-i \eps \tau V(x)} e^{i\frac{\tau}{2}\Delta}\psi^{[n]}(x), \ x \in \overline{\Omega},
\end{equation}
with $ \psi^{[0]}(x) = \psi_0(x)$.

In space, we discretize the Schr\"odinger equation \eqref{eq:linear_1D} by the Fourier pseudospectral method. Let $N$ be an even  positive integer and choose the spatial mesh size $h=(b-a)/N$, then the grid points are  given as
\begin{equation}
x_j := a + jh,\quad j \in \mathcal{T}^0_N=\{j \ \vert \ j = 0, 1, \ldots, N\}.
\end{equation}
Denote $X_N := \{u= (u_0, u_1, \ldots, u_N)^T \in \mathbb{C}^{N+1} \ \vert \ u_0 = u_N\}$ with
the $l^{\infty}$-norm in $X_N$ given as
\begin{equation}
\|u\|_{l^{\infty}} = \max_{0 \le j \le N-1}\vert u_j\vert, \quad u \in X_N.	
\end{equation}
Define $C_{\rm per}(\Omega)=\{u \in C(\overline \Omega) \ \vert \ u(a) = u(b)\}$ and \begin{align*}
& Y_N := \text{span}\left\{e^{i\mu_l(x-a)},\ x \in \overline{\Omega}, \ l \in \mathcal{T}_N\right\},\; \mathcal{T}_N = \left\{l \ \vert \ l = -\frac{N}{2}, \ldots, \frac{N}{2}-1\right\},
\end{align*}
where $\mu_l=\frac{2\pi l}{b-a}$. For any $u(x) \in C_{\rm per}(\Omega)$ and a vector $u \in X_N$, let $P_N: L^2(\Omega) \to Y_N$ be the standard $L^2$-projection operator onto $Y_N$, $I_N : C_{\rm per}(\Omega) \to Y_N$ or $I_N : X_N \to Y_N$ be the trigonometric interpolation operator \cite{ST}, i.e.,

\begin{equation*}
P_N u  = \sum_{l \in \mathcal{T}_N} \widehat{u}_l e^{i\mu_l(x-a)},\quad I_N u  = \sum_{l \in \mathcal{T}_N} \widetilde{u}_l e^{i\mu_l(x-a)},\quad x \in \overline{\Omega},
\end{equation*}
where
\begin{equation*}
\widehat{u}_l = \frac{1}{b- a}\int^{b}_{a} u(x) e^{-i\mu_l (x-a)} dx, \quad \widetilde{u}_l = \frac{1}{N}\sum_{j=0}^{N-1} u_j e^{-i\mu_l (x_j-a)}, \quad l \in \mathcal{T}_N,
\end{equation*}
with $u_j$ interpreted as $u(x_j)$ when involved.

Let $\psi_j^n$ be the numerical approximation of $\psi(x_j,t_n)$ for $j\in \mathcal{T}^0_N$ and $n\ge0$, and denote $\psi^n=(\psi_0^n, \psi_1^n,\ldots, \psi_N^n)^T\in X_N$ as the solution vector. Then, the time-splitting
Fourier pseudospectral (TSFP) method for discretizing the Schr\"odinger equation \eqref{eq:linear_1D}  can be given for $n\ge0$ as
\begin{equation}
\label{eq:TSFP}
\begin{split}
&\psi^{(1)}_j=\sum_{l \in \mathcal{T}_N} e^{-i\frac{\tau\mu^2_l}{2}}\;\widetilde{(\psi^n)}_l\; e^{i\mu_l(x_j-a)}, \\
&\psi^{(2)}_j=e^{-i\eps \tau V(x_j)} \psi^{(1)}_j, \qquad\qquad \qquad   j\in \mathcal{T}^0_N,  \\
&\psi^{n+1}_j=\sum_{l \in \mathcal{T}_N} e^{-i\frac{\tau\mu^2_l}{2}} \; \widetilde{\left(\psi^{(2)}\right)}_l\; e^{i\mu_l(x_j-a)},
\end{split}
\end{equation}
where $\psi^0_j = \psi_0(x_j)$  for $j\in \mathcal{T}^0_N$.
\begin{remark}
The second-order Strang splitting is used for discretizing the Schr\"odinger equation \eqref{eq:linear_1D}. It is straightforward to design the first-order scheme via the Lie splitting and higher order scheme via a higher order splitting method, e.g., the fourth-order compact splitting method or partitioned Runge-Kutta splitting method \cite{BS,MQ,Thal2}.
\end{remark}

\subsection{Local truncation error for the TSFP method}
For proving the (improved) uniform error bounds, we give some results for the local truncation error in this subsection.

We assume the exact solution $\psi(x, t)$ of the Schr\"odinger equation \eqref{eq:linear_1D} up to the time $T_{\varepsilon} = T/\varepsilon$ for any $T>0$ satisfies
\[
{\rm(A)} \quad
\left\|\psi(x, t)\right\|_{L^{\infty}\left([0, T_{\varepsilon}]; H^{m}_{\rm per}\right)} \lesssim 1, \quad  \left\|\partial_t \psi(x, t)\right\|_{L^{\infty}\left([0, T_{\varepsilon}]; H^{m-2}_{\rm per}\right)} \lesssim 1,
\]
and the potential satisfies
\[
{\rm(B)} \quad
V(x) \in  H^{m^{\ast}}_{\rm per}, \quad m^{\ast} = \max\{m, 5\},
\]
where $m$ describes the regularity of the exact solution. Here, $H_{\rm per}^m(\Omega)=\{\phi \in H^m(\Omega) \vert \partial_x^k\phi(a)=\partial_x^k\phi(b),\,k=0,1,\ldots,m-1\}$, with the equivalent $H^m$-norm on $H_{\rm per}^m(\Omega)$ given as
$\|\phi\|_{H^m}=\left(\sum\limits_{l\in\mathbb{Z}}(1+\mu_l^{2})^m \vert \widehat{\phi}_l\vert^2\right)^{1/2}$.
In the rest of this paper, we may write $\psi(t) = \psi(x, t)$, i.e. omit the spatial variable, when there is no confusion.  The following estimates of the local truncation error for the semi-discretization \eqref{eq:Strang} hold.
\begin{lemma}
\label{lemma:local-linear}
Under assumptions (A) and (B) with $m\ge3$, for $0 < \eps \leq 1$, the local truncation error of the  TSFP  \eqref{eq:TSFP} for the Schr\"odinger equation with $O(\eps)$-potential at time $t_n$ can be written as ($0\leq n\leq T_\eps/\tau-1$)
\begin{equation}\label{eq:localerr}
\mathcal{E}^{n}(x) := P_N\mathcal{S}_{\tau}(P_N\psi(t_n)) -P_N \psi(t_{n+1}) = P_N\mathcal{F}(P_N\psi(t_n))+ R_{n},
\end{equation}
where
\begin{equation}\label{eq:Fpsi}
\mathcal{F}(P_N\psi(t_n))  = -i \varepsilon\tau f^n\left(\frac{\tau}{2}\right) + i\varepsilon\int^{\tau}_0 f^n(s)\, ds,
\end{equation}
with $f^n(s)=e^{i(\tau-s)\Delta} V e^{is\Delta}P_N\psi(t_n)$, and
the following error estimates hold
\begin{equation}
\|\mathcal{F}(P_N\psi(t_n))\|_{H^k}\lesssim \eps \tau^3,\quad \|R_{n}\|_{H^k} \lesssim \varepsilon^2\tau^3+\varepsilon\tau h^{m-k},\quad k=0,1.
\label{eq:local_linear}
\end{equation}
In addition, $\mathcal{E}^{n}(x),  R_n(x)\in Y_N$ and the $L^2$-estimates in \eqref{eq:local_linear} hold for $m\ge2$.
\end{lemma}
\noindent
\begin{proof}
The proof is standard following \cite{JL00,LU}, and we sketch the procedure to emphasize the effects of spatial discretization and the
parameter $\varepsilon$.
 By the Taylor expansion for $e^{-i\varepsilon\tau V}$, we have
\begin{align*}
P_N\left(\mathcal{S}_{\tau}(P_N\psi(t_n))\right) = &\ e^{i\tau\Delta}P_N\psi(t_n)- i\varepsilon\tau P_N\left(e^{i\frac{\tau}{2}\Delta} V  e^{i\frac{\tau}{2}\Delta} P_N\psi(t_n)\right)  \\
&\   - \varepsilon^2\tau^2P_N\left(\int^1_0 (1-\theta) e^{i\frac{\tau}{2}\Delta}e^{-i\varepsilon \theta\tau V}  V^2 e^{i\frac{\tau}{2}\Delta}P_N\psi(t_n) d\theta\right).
\end{align*}
On the other hand, by repeatedly using the Duhamel's principle, we can write
\begin{align}
P_N\psi(t_{n+1}) = & \ P_N\left(e^{i\tau\Delta}	\psi(t_n)\right) - i\varepsilon P_N\left(\int^{\tau}_0 e^{i(\tau-s)\Delta} V  e^{is\Delta}\psi(t_n) d s\right) \nonumber\\
&\ - \varepsilon^2 P_N\left(\int^{\tau}_0 \int^s_0 e^{i(\tau-s)\Delta}V e^{i(s-w)\Delta} V\psi(t_n + w) d w d s\right).\nonumber
\end{align}
Recalling assumptions (A) and (B), applying Fourier projections, we  have
\begin{align*}
P_N\psi(t_{n+1})	 = &\ e^{i\tau\Delta}	P_N\psi(t_n) - i\varepsilon \int^{\tau}_0 P_N\left(e^{i(\tau-s)\Delta} V  e^{is\Delta}P_N\psi(t_n) \right)d s \\
& \ - \varepsilon^2 \int^{\tau}_0 \int^s_0 P_N\left(e^{i(\tau-s)\Delta}V e^{i(s-w)\Delta} VP_N\psi(t_n + w) \right)d w d s
-r_h^n,
\end{align*}
with $\|r_h^n(x)\|_{L^2}\lesssim \varepsilon \tau h^m$ and $\|r_h^n(x)\|_{H^1}\lesssim \varepsilon\tau h^{m-1}$.
Introducing $f^n(s)$ as in Lemma \ref{lemma:local-linear} and
\begin{equation*}
B^n\left(s, w\right) =P_N\left( e^{i(\tau-s)\Delta}V e^{i(s-w)\Delta}Ve^{iw\Delta}P_N\psi(t_n)\right), \quad 	
\end{equation*}
the local truncation error can be written as \cite{LU}
\begin{align*}
\mathcal{E}^{n}= & \ P_N\mathcal{F}(P_N\psi(t_n)) - \frac{\varepsilon^2\tau^2}{2} B^n\left(\frac{\tau}{2}, \frac{\tau}{2}\right) + \varepsilon^2 \int^{\tau}_0\int^s_0 B^n\left(s, w\right) d w d s \\
&\ + \varepsilon^2 r_1^n + \varepsilon^2 r_2^n+r_h^n,
\end{align*}
where $\mathcal{F}(P_N\psi(t_n))$ is given in \eqref{eq:Fpsi} and
\begin{align*}
& r_1^n = -\tau^2 \int^1_0 (1-\theta) P_N\left(e^{i\frac{\tau}{2}\Delta}  (e^{-i\varepsilon \theta\tau V}-1) V^2 e^{i\frac{\tau}{2}\Delta}P_N\psi(t_n)\right) d\theta,   \\
& r_2^n =  \int^{\tau}_0 \int^s_0\left(P_N\left( e^{i(\tau-s)\Delta}V e^{i(s-w)\Delta} VP_N\psi(t_n + w) \right)-   B^n\left(s, w\right) \right)d w d s.	
\end{align*}
Since $e^{i\tau\Delta}$ preserves the $H^s$-norm and $\|(e^{-i\varepsilon \theta\tau V}-1) V^2 \|_{H^1}\lesssim \varepsilon\tau\theta\|V\|_{H^1}^3$, we have
\begin{align*}
\left\|r_1^n\right\|_{H^1} &  \lesssim \varepsilon\tau^3\|V\|_{H^1}^3\|\psi(t_n)\|_{H^1} \lesssim \varepsilon\tau^3.
\end{align*}
The following estimates are standard  (c.f. \cite{LU}),
\begin{align*}
&\left\|r_2^n\right\|_{H^1} 
 \lesssim \varepsilon\tau^3\|V\|_{H^1}^2 \left\|V \psi(\cdot)\right\|_{L^{\infty}([0, \tau]; H^1)}\lesssim  \varepsilon\tau^3,\\
& \left\|-\frac{\tau^2}{2} B\left(\frac{\tau}{2}, \frac{\tau}{2}\right) + \int^{\tau}_0\int^s_0  B\left(s, w\right) d w d s \right\|_{H^1}
\lesssim \tau^3 \|V\|_{H^3}^2\|\psi(t_n)\|_{H^3}\lesssim\tau^3.
\end{align*}
Finally, for the major part of the local truncation error can be estimated by the midpoint quadrature rule as
\cite{LU}
\begin{equation}
\|\mathcal{F}(P_N\psi(t_n))\|_{H^1}\lesssim \varepsilon\tau^3\|[\Delta, [\Delta, V]]P_N\psi(t_n)\|_{H^1} \lesssim \|V\|_{H^5}{\|\psi(t_n)\|_{H^3}},
\end{equation}
where $[\Delta, [\Delta, V]]$ is the double commutator. Thus,  by setting
\begin{equation}
 R_{n}=- \frac{\varepsilon^2\tau^2}{2} B^n\left(\frac{\tau}{2}, \frac{\tau}{2}\right) + \varepsilon^2 \int^{\tau}_0\int^s_0 B^n\left(s, w\right) d w d s + \varepsilon^2 r_1^n + \varepsilon^2 r_2^n+r_h^n,
\end{equation}
we obtain the estimates in Lemma \ref{lemma:local-linear}.
\end{proof}

\subsection{Uniform error bounds in $L^2$-norm}
In this subsection, we adopt the unitarity of the numerical solution flow in $L^2(\Omega)$ to establish the uniform error bound  in $L^2$-norm with linear growth in $t$ up to the time $t\leq T_\eps=\frac{T}{\varepsilon}$. We remark that the uniform estimates are standard, while the linear growth of the error only holds in the $L^2$-norm. The reason is that TSFP \eqref{eq:TSFP} only preserves the $L^2$-norm.
\begin{theorem}
\label{thm:eb_linear}
Let $\psi^n$ be the numerical approximation obtained from the TSFP \eqref{eq:TSFP}. Under  assumptions (A) and (B) with $m \ge 2$, for any $0 < \varepsilon \leq 1$, we have
\begin{equation}
\left\|\psi(x, t_n) - I_N\psi^n\right\|_{L^2} \leq (C_0 + C_1\eps t_n) \tilde{C}(T)\left(h^m +  \tau^2\right), \quad  0 \leq n \leq \frac{T/\varepsilon}{\tau}, \\
\label{eq:error_l}
\end{equation}
where $C_0$ and $C_1$ are two positive constants independent of  $h$, $\tau$, $n$, $\eps$ and $T$, $\tilde{C}(T)$ depends on $\|\psi\|_{L^\infty([0,T];H^m)}$ and $\|V\|_{H^{m^*}}$.
\end{theorem}	
\noindent
\begin{proof}
 Noticing that
\begin{equation}
I_N \psi^n - \psi(t_n) =I_N \psi^n - P_N(\psi(t_n)) +   P_N(\psi(t_n)) - \psi(t_n),
\end{equation}
under assumptions (A) and (B), we get from the standard Fourier projection properties \cite{ST}
\begin{equation}\label{eq:dec}
\left\|I_N \psi^n - \psi(t_n)\right\|_{L^2} \leq \left\|I_N \psi^n - P_N(\psi(t_n))\right\|_{L^2} + C_2 h^{m},\ 0 \leq n \leq \frac{T/\eps}{\tau}.
\end{equation}
Thus, it suffice to consider the error function $e^n\in Y_N$ at $t_n$ as
\begin{equation}\label{eq:err-def}
e^n:=e^n(x) = I_N\psi^n -P_N \psi(t_n), \quad 0 \leq n \leq \frac{T/\eps}{\tau},
\end{equation}
 and $\left\|e^0\right\|_{L^2} \leq C_3 h^m$ implied by the standard projection and interpolation results. From the local error \eqref{eq:localerr} in Lemma \ref{lemma:local-linear}, we have the error equation for $e^{n}$ ($0 \leq n \leq \frac{T/\eps}{\tau}-1$),
 \begin{equation}\label{eq:errorg}
 e^{n+1}=I_N\psi^{n+1}-P_N\psi(t_{n+1})=I_N\psi^{n+1}- P_N\mathcal{S}_{\tau}(P_N\psi(t_n))  + \mathcal{E}^{n}.
 \end{equation}
Noticing the fully discrete scheme \eqref{eq:TSFP} and $\mathcal{S}_\tau$ \eqref{eq:Strang}, i.e.
\begin{align*}
&I_N \psi^{n+1} = e^{i\frac{\tau}{2}\Delta}(I_N\psi^{(2)}),\, I_N(\psi^{(2)}) =  I_N(e^{-i\eps \tau V(x)}\psi^{(1)} ),\, I_N\psi^{(1)}=e^{i\frac{\tau}{2}\Delta} I_N\psi^{n},\\
&P_N(\mathcal{S}_\tau(\psi(t_n)))=e^{i\frac{\tau}{2}\Delta}(P_N\psi^{\langle2\rangle}),\,
\psi^{\langle2\rangle}=e^{-i\eps \tau V(x)} \psi^{\langle1\rangle},\,\psi^{\langle1\rangle}=e^{i\frac{\tau}{2}\Delta} P_N\psi(t_n),
\end{align*}
in view of the facts that $I_N$ and $P_N$ are identical on $Y_N$ and $e^{i\tau\Delta/2}$ preserves the $H^k$-norm ($k\ge0$), using Taylor expansion $e^{-i\eps\tau V(x)}=1-i\eps\tau V(x)\int_0^1e^{-i\eps\theta\tau V(x)}\,d\theta$ and assumptions (A) and (B),
we have
\begin{align}
&\|I_N\psi^{n+1}- P_N\mathcal{S}_{\tau}(P_N\psi(t_n))\|_{L^2}=\|I_N\psi^{(2)}-P_N\psi^{\langle2\rangle}\|_{L^2},\label{eq:err-1} \\
&\|P_N\psi^{\langle2\rangle}-I_N\psi^{\langle2\rangle}\|_{L^2}=\left\|\eps\tau(P_N-I_N) \left(V(x)\int_0^1e^{-i\eps\theta\tau V(x)}\,d\theta\psi^{\langle1\rangle}\right)\right\|_{L^2}\nonumber\\
&\qquad\qquad\qquad\quad\qquad\;\, \leq C_4\eps\tau h^m,\label{eq:err-2}
\end{align}
where $C_4$ is obtained from Fourier interpolation and projection properties together with
$\left\| \left(V(x)\int_0^1e^{-i\eps\theta\tau V(x)}\,d\theta\psi^{\langle1\rangle}\right)\right\|_{H^m}\leq C(\|V\|_{H^m})\|\psi(t_n)\|_{H^m}$. In addition, by direct computation and Parseval's identity, we can derive
\begin{align}
\|I_N\psi^{(2)}-I_N\psi^{\langle2\rangle}\|_{L^2}&=\sqrt{h\sum\limits_{j=0}^{N-1}\vert \psi^{(2)}_j-\psi^{\langle2\rangle}(x_j)\vert^2}=\sqrt{h\sum\limits_{j=0}^{N-1}\vert\psi^{(1)}_j-\psi^{\langle1\rangle}(x_j)\vert^2}\nonumber\\
&=\|I_N\psi^{(1)}-I_N\psi^{\langle1\rangle}\|_{L^2}=\|I_N\psi^n-P_N\psi(t_n)\|_{L^2}\nonumber\\
&=\|e^n\|_{L^2}.\label{eq:err-3}
\end{align}
Taking the $L^2$-norm on both sides of \eqref{eq:errorg} and combining \eqref{eq:err-1},\eqref{eq:err-2} and \eqref{eq:err-3} together,  in view of Lemma \ref{lemma:local-linear}, we obtain for $0\leq n\leq\frac{T/\eps}{\tau}-1$,
\begin{align}
\|e^{n+1}\|_{L^2} &\leq \|\mathcal{E}^{n}\|_{L^2}+\|I_N\psi^{(2)}-P_N\psi^{\langle2\rangle}\|_{L^2}\nonumber\\
& \leq  \|\mathcal{E}^{n}\|_{L^2}+\|I_N\psi^{(2)}-I_N\psi^{\langle2\rangle}\|_{L^2}+\|P_N\psi^{\langle2\rangle}-I_N\psi^{\langle2\rangle}\|_{L^2}\nonumber\\
& \leq  \|e^n\|_{L^2} + C_5\left(\varepsilon\tau h^m+\eps\tau^3\right).\label{eq:err-4}
\end{align}
Thus, the following estimates hold
\begin{equation}
\|e^{n+1}\|_{L^2}\leq C_5\eps t_{n+1}( h^m+\tau^2)+C_3h^m,\quad 0\leq n\leq\frac{T/\eps}{\tau}-1,
\end{equation}
and the conclusion of Theorem \ref{thm:eb_linear} by taking $C_0 = C_2+C_3$ and $C_1 = C_5$ in view of \eqref{eq:dec}. It is easy to verify all the constants appearing in the proof only depend on $V$ and $\psi$.
\end{proof}

\begin{remark} Regarding the estimate \eqref{eq:error_l} in Theorem \ref{thm:eb_linear},  $\tilde{C}(T)$ comes from the local truncation error, which depends on the growth of the Sobolev norm w.r.t $T$ in the assumption (A).  Based on previous analytical results, $\tilde{C}(T)$  usually has a polynomial growth in $T$ \cite{Bour2}, and $\tilde{C}(T)$ could be uniformly bounded w.r.t. $T$ for certain type of potential function $V(x)$ \cite{Wang}.
\end{remark}
\begin{remark}
According to Theorem \ref{thm:eb_linear}, the uniform error bound for the TSFP method in $L^2$-norm at time $t_\eps=t/\eps$ for the Schr\"odinger equation linearly grows with respect to $t$, and the results can be generalized to other splitting methods.
In fact, given an accuracy bound $\delta_0 > 0$, the time (for simplicity, assume $\eps=1$ here) for the second-order splitting method to violate the accuracy requirement  $\delta_0$ is $O(\delta_0/\tau^2)$. For the first-order and fourth-order splitting methods, the time is $O(\delta_0/\tau)$ and $O(\delta_0/\tau^4)$, respectively. In other words, higher order splitting method performs much better in the long-time simulations not only regarding the higher accuracy but also longer simulation time to produce accurate solutions. For the
 $L^2$-estimates in Theorem \ref{thm:eb_linear}, the regularity requirements on the potential $V(x)$  can be weakened. In addition, extensions to 2D/3D are straightforward.
\end{remark}
\begin{remark}
By the similar procedure (or formally letting $h\to0^+$), we could establish uniform error bounds for the semi-discretization.
Let $\psi^{[n]}$ be the numerical approximation obtained from the Strang splitting \eqref{eq:Strang}. Under  assumptions (A) and (B) with $m \ge 3$, for any $0 < \varepsilon \leq 1$, we have
\begin{equation}
\|\psi(x, t_n) - \psi^{[n]}\|_{L^2} \leq  C_0 \eps t_n\tau^2, \quad  0 \leq n \leq \frac{T/\varepsilon}{\tau}, \\
\label{eq:error_l2}
\end{equation}
where $C_0$ is a positive constant independent of $\tau$, $n$ and $\eps$. Such linear growth of the error constant w.r.t $t_n$ in \eqref{eq:error_l2} has been  previously reported in \cite{JL00}.
\end{remark}

\subsection{Improved uniform error bounds in $H^1$-norm}
In this subsection, we show  improved uniform error bounds in $H^1$-norm for the Schr\"odinger equation with $O(\eps)$-potential up to the time $T_\eps = T/\eps$ under assumptions (A) and (B) with $m \ge 3$, where
we will work with $H^1$-estimates for the nonlinear case also to control the nonlinearity in 1D. It is worth noticing that in higher dimensions (2D/3D), $H^2$-estimates would be enough.
 The improved estimates rely on the cancellation phenomenon of non-resonant oscillating frequencies, for which we shall require the time step size $\tau$ satisfy certain non-resonant conditions. In the fully discrete case, for the Fourier modes $|l|\leq \lceil\frac{1}{\tau_0}\rceil$ ($\tau_0\in(0,1)$, $\lceil\cdot\rceil$ is the  ceiling function), we impose the Diophantine type condition \cite{DG,SZ}: there exists a constant  $C_0>0$ such that
\begin{equation}\label{eq:nonres}
\left|1-e^{i\tau \mu_1^2 K}\right|\ge \frac{C_0\tau^{\nu_1}}{(\mu_1^2 |K|)^{\nu_2}},\quad 0<|K|\leq K_0=\lceil1/\tau_0\rceil^2,\quad K\in\mathbb{Z},
\end{equation}
where $\nu_1\in[0,1]$, $\nu_2\ge -1$, and the bound $|K|\leq \lceil1/\tau_0\rceil^2$ corresponds to the interaction between potential $V(x)$ and the solution $\psi(x,t)$. In particular, we consider the following cases of time step sizes: for a given constant $\alpha\in(0,1)$, the time step size $\tau$ satisfies
\begin{equation}\label{eq:tauc1}
\tau\in\left(0, \alpha\frac{2\pi}{\mu_1^2(1+\tau_0)^2}\tau_0^2\right),\end{equation}
or the Diophantine type step condition \cite{SZ} 
\begin{equation}\label{eq:tauc2}
\tau\in I_{\nu_3,\alpha,\tau_0}=\left\{\tau>0: \left|\tau-\frac{2l\pi}{\mu_1^2K}\right|\ge\frac{\lambda}{|\mu_1^2K|^{2+\nu_3}},  K,l\in\mathbb{Z},\; 0<|K|\leq K_0,\;
0\leq l\right\},
\end{equation}
where $\lambda=\frac{\alpha\pi\mu_1^{2+2\nu_3}}{4\sum\limits_{k=1}
^\infty1/k^{1+\nu_3}}$ and $\nu_3>0$. \eqref{eq:tauc2} is adapted from a general form in \cite{SZ}, and it is direct to observe that 
\begin{equation*}
[0,2\pi/\mu_1^2]\backslash I_{\nu_3,\alpha,\tau_0}=\bigcup_{0\leq l\leq K, 1\leq K\leq K_0}\left\{\tau\in[0,2\pi/\mu_1^2]:\left|\tau-\frac{2l\pi}{\mu_1^2K}\right|<\frac{\lambda}{|\mu_1^2K|^{2+\nu_3}}\right\},
\end{equation*}
where the Lebesgue measure of RHS is bounded by $\alpha\pi/(2\mu_1^2)$ and $I_{\nu_3,\alpha,\tau_0}\cap [0,2\pi/\mu_1^2]$ has measure greater than $3\pi/2\mu_1^2$. Moreover, if $\tau\in I_{\nu_3,\alpha,\tau_0}$, $\tau+2k\pi\in I_{\nu_3,\alpha,\tau_0}$ ($k\ge0,k\in\mathbb{Z}$) and large time step sizes are admissible in \eqref{eq:tauc2}.

Now we can verify that \eqref{eq:tauc1} fulfills \eqref{eq:nonres} with $\nu_1=1,\nu_2=-1,C_0=\frac{\sin(\alpha\pi)}{\pi\alpha}$, while \eqref{eq:tauc2} fulfills \eqref{eq:nonres} with $\nu_1=0,\nu_2=1+\nu_3,C_0=C_{\nu_3}\alpha$ ($C_{\nu_3}$ a constant depending on $\nu_3$, see also \eqref{eq:cfl2}).
\eqref{eq:tauc1} corresponds to the typical choice of time step size $\tau$ allowing $\tau\to0^+$ and \eqref{eq:tauc2} allows $\varepsilon$ dependent large time step size which is well suited for the long time  dynamics of Schr\"odinger equation \eqref{eq:linear}. We remark here that similar non-resonance condition was used for establishing uniform error bounds of time-splitting methods for the (nonlinear) Dirac equation in the nonrelativistic regime \cite{BCY20,BCY21}.
We refer to Remark \ref{rmk:tau0} for more discussions on the non-resonance condition \eqref{eq:nonres}.

Under the non-resonance condition \eqref{eq:nonres}, we have the following improved estimates.

 \begin{theorem}
\label{thm:im_linear}
Let $\psi^n$ be the numerical approximation obtained from the TSFP \eqref{eq:TSFP}. Under the assumptions (A) and (B) with $m\ge3$, for any $\varepsilon \in (0, 1]$ and a fixed $\tau_0\in(0,1)$, when $\tau$ satisfies \eqref{eq:tauc1} or \eqref{eq:tauc2} ($m\ge 5+2\nu_3$ for \eqref{eq:tauc2} case), we have the estimates

\begin{equation}
\left\|\psi(x, t_n) - I_N\psi^n\right\|_{H^1} \lesssim h^{m-1} + \varepsilon \tau^2 +\tau_0^{m-1}, \quad  0 \leq n \leq \frac{T/\varepsilon}{\tau}.
\label{eq:error_2}
\end{equation}
In particular,  if the exact solution is smooth, i.e.  $\psi(x, t) \in H^{\infty}_{\rm per}$,  the  $\tau_0^{m-1}$ part error would decrease exponentially in terms of $\tau_0$ and  can be ignored in practical
computation when $\tau_0$ is taken as $\tau_0^{\rm cr}$ (small but fixed,  only depends on $\psi$ and the logarithm of the machine precision), thus the improved error bounds for sufficiently small $\tau$ could be stated as
\begin{equation}
\left\|\psi(x, t_n) - I_N\psi^n\right\|_{H^1} \lesssim h^{m-1} + \varepsilon \tau^2, \quad  0 \leq n \leq \frac{T/\varepsilon}{\tau}.
\label{eq:error_3}
\end{equation}
\end{theorem}
 \begin{remark}
Before the presentation of the proof, some observations are marked.

1.  First, $1/\tau_0>1$  serves as a cut-off mode, i.e. the modes $\vert l \vert >1/\tau_0$ are treated by Fourier projection, and the modes $\vert l \vert  \leq\lceil\frac{1}{\tau_0}\rceil$ will be treated by the RCO technique for non-resonant $\tau$ in \eqref{eq:tauc1}-\eqref{eq:tauc2}. 

 2. For  \eqref{eq:tauc1}, the requirement is that the step size $\tau$ resolves the largest oscillatory frequency of the free Schr\"odinger operator below the cut-off modes as $\vert \mu_l\vert^2=\frac{4l^2\pi^2}{(b-a)^2}\sim \frac{4\pi^2}{(b-a)^2\tau_0^2}$, i.e. $\tau\vert \mu_l\vert^2<2\pi$. In turn, the error constant in front of $\varepsilon\tau^2$ depend on the parameter $\alpha\in(0,1)$ (scales like $\frac{\alpha\pi}{\sin(\alpha\pi)}$). Thus, the introduced parameter $\tau_0$ can be either fixed, or any other choices satisfying the condition $\tau <  \alpha\frac{2\pi}{\mu_1^2(1+\tau_0)^2}\tau_0^2$, e.g. $\tau_0=\frac{\sqrt{2}\mu_1}{\sqrt{\pi\alpha}}\sqrt{\tau}$.  $\tau_0$ serves as a Fourier projection parameter, similar to the role of the spatial mesh size $h$.
Alternatively, we can also choose $\tau_0 = 2/N$ such that the last term in \eqref{eq:error_2} could be controlled by the first term, and the requirement \eqref{eq:tauc1} on $\tau$ becomes a CFL type condition $\tau\lesssim h^2$.

3. For the larger non-resonance step size $\tau$ in \eqref{eq:tauc2}, Theorem \ref{thm:im_linear} implies that for a given accuracy  $\delta_0$,  $\tau$ can be chosen as $O(\sqrt{\delta_0}/\sqrt{\varepsilon})$ large, which is particularly superior for $\varepsilon\ll1$. We notice that \eqref{eq:tauc2} require higher regularity for deriving the improved error bounds.

4. For general non-resonance  time step sizes satisfying \eqref{eq:nonres}, the improved error estimates \eqref{eq:error_2} hold by the similar arguments with slightly different regularity assumptions on the exact solution $\psi(x,t)$. Moreover, the error constant in front of the $\varepsilon\tau^2$ term depends on $C_0$ in \eqref{eq:nonres} as $\sim 1/C_0$.

5. For 2D/3D extensions, the improved error bounds can be directly established for the step sizes in \eqref{eq:tauc1} and non-resonance step sizes as  \eqref{eq:tauc2} in higher dimensions. See Remark \ref{remark:hd} for more details.
\end{remark}
\begin{proof}
Following the proof of Theorem \ref{thm:eb_linear}, we only need to estimate the error $e^n$ in \eqref{eq:err-def} for $0\leq n\leq\frac{T/\eps}{\tau}$. First, using the fact that $P_N=I_N$ when it is restricted on $Y_N$, we can write for $0 \leq n \leq \frac{T/\eps}{\tau}$,
\begin{align}
&I_N\psi^{n+1}- P_N\mathcal{S}_{\tau}(P_N\psi(t_n))=e^{i\tau \Delta}\left(I_N\psi^{n}-P_N\psi(t_n)\right)+Q^n(x),\label{eq:error-g}
\end{align}
where $Q^n(x)\in Y_N$ is given by
\begin{align}
Q^n(x) &=-i\eps \tau e^{i\frac{\tau}{2}\Delta}\left(I_N\left(V(x)\int_0^1e^{-i\eps\theta\tau V(x)}\,d\theta \psi^{(1)}\right)\right)\nonumber\\
&\quad +i\eps\tau  e^{i\frac{\tau}{2}\Delta}\left(P_N\left(V(x)\int_0^1e^{-i\eps\theta\tau V(x)}\,d\theta \psi^{\langle1\rangle}\right)\right).\label{eq:Q-def}
\end{align}
Using Parseval's identity and finite difference operator (c.f. \cite{BC,BC14}), by similar estimates \eqref{eq:err-1}, \eqref{eq:err-2} and \eqref{eq:err-3}  for the $L^2$-norm case,  we can control $Q^n$ as
\begin{equation}\label{eq:Qbd}
\|Q^n\|_{H^1}\lesssim\eps\tau\left(h^{m-1}+\|e^n\|_{H^1}\right),\quad 0\leq n\leq\frac{T/\eps}{\tau}-1.
\end{equation}
From \eqref{eq:error-g} and \eqref{eq:errorg}, we could derive that for $0 \leq n \leq \frac{T/\eps}{\tau}-1$,
\begin{equation}
e^{n+1}
 =  \ e^{i\tau\Delta}e^n + Q^n(x) + \mathcal{E}^{n},
 \label{eq:etg}
\end{equation}
which implies
\begin{equation}
e^{n+1} = e^{i(n+1)\tau\Delta}e^0 + \sum^{n}_{k = 0}e^{i(n-k)\tau \Delta}\left(Q^k(x) + \mathcal{E}^{k}\right).
\end{equation}
{\bf{Step 1.}} (Identifying the leading error term) Using the local truncation error representation \eqref{eq:local_linear} in Lemma \ref{lemma:local-linear}, we have
\begin{equation}
 \sum^{n}_{k = 0}e^{i(n-k)\tau \Delta}\mathcal{E}^{k}  = \sum^{n}_{k = 0}e^{i(n-k)\tau \Delta}(P_N\mathcal{F}(P_N\psi(t_k)) + R_{k}), 	
\end{equation}
and
\begin{align}
& \left\|\sum^{n}_{k = 0}e^{i(n-k)\tau \Delta}R_{k}\right\|_{H^1} \lesssim (n+1) \left(\eps^2 \tau^3+\eps\tau h^{m-1}\right) \lesssim T \eps \tau^2+T h^{m-1}, \\
& \left\|\sum^{n}_{k = 0}e^{i(n-k)\tau \Delta} Q^{k}(x)\right\|_{H^1}  \lesssim \eps\tau\sum^{n}_{k=0} \left\|e^k\right\|_{H^1} + h^{m-1}.
\end{align}
Combining above estimates and $\|e^0\|_{H^1} \lesssim h^{m-1}$, we obtain for $0 \leq n \leq \frac{T/\eps}{\tau}-1$,
\begin{align}
\|e^{n+1}\|_{H^1} \lesssim & \ h^{m-1} + \eps \tau^2 + \eps\tau\sum^{n}_{k=0} \|e^k\|_{H^1}\nn \\
 & \ + \left\|\sum^{n}_{k = 0}e^{i(n-k)\tau \Delta}P_N\mathcal{F}(P_N\psi(t_k))\right\|_{H^1}.
\label{eq:S1eg}
\end{align}
Recalling Lemma \ref{lemma:local-linear},  we have $\|\mathcal{F}(P_N\psi(t_k))\|_{H^1}\lesssim \tau^3$, which implies $\|e^{n+1}\|_{H^1}\lesssim \tau^2+h^{m-1}$.
Thus, to prove the improved error estimates, we need analyze  the last term in \eqref{eq:S1eg} carefully, i.e., treat the sum $\sum^{n}_{k = 0}e^{i(n-k)\tau \Delta}P_N\mathcal{F}(P_N\psi(t_k))$ in a proper way.
To gain an order of $O(\eps)$ from the sum, we shall introduce the {\bf regularity compensated oscillation} (RCO) technique. From \eqref{eq:linear}, we find $\partial_t\psi(x,t)-i\Delta\psi(x,t)=O(\eps)$, and it is natural to introduce the `twisted variable' as
 \begin{equation}
 \phi(x, t) = e^{-it\Delta}\psi(x, t), \quad t \geq 0,
 \label{eq:twist1}
 \end{equation}
 and  $\phi(t) := \phi(x, t)$ satisfies the equation
 \begin{equation}\label{eq:twist}
 i\partial_t\phi(x,t)=\eps e^{-it\Delta}\left(V(x)e^{it\Delta}\phi(x,t)\right),\quad t > 0.
 \end{equation}
 It is direct to see that $\phi(x,t)$ enjoys the same $H^k$ ($k\ge0$) bounds as $\psi(x, t)$, while
\begin{equation}
\left\|\partial_t\phi(t)\right\|_{H^{m}} \lesssim \eps, \quad 0\leq t\leq T/\eps.\label{eq:S1twist}
\end{equation}
The RCO approach would then perform a summation-by-parts procedure in the $\sum^{n}_{k = 0}e^{i(n-k)\tau \Delta}P_N\mathcal{F}(P_N\psi(t_k))$ to force $\partial_t\phi(t)$ appear with a gain of order $O(\eps)$, where $\tau$ is small to control the accumulation of the phase (frequency) of the type $e^{i(n-k)\tau \Delta}$. Since the number $N$ of the spatial grid points   could be very large, we shall introduce a cut-off parameter $\tau_0\in(0,1)$, where the high frequency modes ($|l|>\frac{1}{\tau_0}$) will be controlled by the smoothness of the exact solution and the Fourier projections, and the low frequency modes ($|l|\leq\frac{1}{\tau_0}$) will be dealt with the RCO technique.

{\bf Cut-off parameter}.  Choose $\tau_0\in(0,1)$, and let $N_0=2\lceil1/\tau_0\rceil\in\mathbb{Z}^+$  with $1/\tau_0\leq N_0/2<1+1/\tau_0$,  then only those Fourier modes with $-\frac{N_0}{2}\leq l\leq \frac{N_0}{2}-1$ in $\mathcal{F}(P_N\psi(t_k))$ would be considered.  Based on the Fourier projections and the assumption (A), we have $\|P_{N_0}\psi(x,t)-P_N\psi(x,t)\|_{L^\infty([0,T/\eps];H^1)}\lesssim h^{m-1}+N_0^{1-m}\lesssim h^{m-1}+\tau_0^{m-1}$ and for $0\leq n\leq \frac{T/\eps}{\tau}-1$,
\begin{equation}\label{eq:cut-off}
\left\|P_{N_0}\mathcal{F}(P_{N_0}\psi(t_n))-P_N\mathcal{F}(P_N\psi(t_n))\right\|_{H^1}
\lesssim \eps\tau (h^{m-1}+\tau_0^{m-1}).
\end{equation}
 Indeed, since $P_N\psi(t_k)\in Y_N$, we could actually assume the choice of $\tau_0$ such that  $N_0\leq N$, but here we work without this condition for the convenience of extension to the semi-discretization-in-time case.

Based on \eqref{eq:S1eg}, \eqref{eq:twist1} and  \eqref{eq:cut-off}, recalling the unitary properties of $e^{it\Delta}$, we find for $0 \leq n \leq \frac{T/\varepsilon}{\tau}- 1$,
\begin{align}
\left\|e^{n+1}\right\|_{H^1} & \lesssim  h^{m-1} + \tau_0^{m-1}+ \eps \tau^2 + \eps\tau\sum^n_{k=0} \left\|e^k\right\|_{H^1}+\|\mathcal{R}^n\|_{H^1}, \label{eq:final1} \\
\mathcal{R}^n(x)& = \sum\limits_{k=0}^n e^{-i(k+1)\tau\Delta} P_{N_0} \mathcal{F}(e^{i t_k\Delta}(P_{N_0}\phi(t_k))) .
\label{eq:Rn}
\end{align}

\noindent
{\bf{Step 2.}} (Analysis  via RCO)  Let $\phi(t) = \sum_{l\in\mathbb{Z}}\widehat{\phi}_l(t)e^{i\mu_l(x-a)}$ ($t\ge0$), and we have $P_{N_0}\phi(t) = \sum_{l\in\mathcal{T}_{N_0}}\widehat{\phi}_l(t)e^{i\mu_l(x-a)}$, where $\widehat{\phi}_l(t)$ is the $l$-th Fourier coefficient of $\phi(x,t)$.
For $l\in\mathcal{T}_{N_0}$, introduce the multi-index set $\mathcal{I}_l^{N_0}$ associated with $l$ as
\begin{equation}
\mathcal{I}_l^{N_0} = \left\{(l_1,l_2) \ \vert \ l_1+l_2 = l,\ l_1\in\mathbb{Z},\, l_2 \in\mathcal{T}_{N_0}\right\}.
\end{equation}
According to the definition of $\mathcal{F}$ in Lemma \ref{lemma:local-linear}, we have the expansion
\begin{equation*}
 e^{-i(k+1)\tau\Delta}P_{N_0}\left(e^{i(\tau-s)\Delta} V e^{is\Delta}P_{N_0}\psi(t_n)\right)=\sum\limits_{l\in\mathcal{T}_{N_0}}\sum\limits_{(l_1,l_2)\in\mathcal{I}_l^{N_0}}
\mathcal{G}_{k,l,l_1,l_2}(s)e^{i\mu_l(x-a)},
\end{equation*}
where $\mathcal{G}_{k,l,l_1,l_2}(s)$  ($l,l_2\in\mathcal{T}_{N_0}$) is a function of $s$ as
\begin{equation}
\mathcal{G}_{k,l,l_1,l_2}(s) =
e^{i(t_k+s )\delta_{l, l_2}}\widehat{V}_{l_1}\widehat{\phi}_{l_2}(t_k), \quad \delta_{l, l_2} = \delta_l - \delta_{l_2},\quad \delta_l = \mu^2_l. \label{eq:mGdef}\\
\end{equation}
Then, the remainder term $\mathcal{R}^n(x)$ in \eqref{eq:final1} reads
\begin{equation}
\label{eq:remainder-dec1}
\mathcal{R}^n(x) = i \eps \sum\limits_{k=0}^n
\sum\limits_{l\in\mathcal{T}_{N_0}}\sum\limits_{(l_1,l_2)\in\mathcal{I}_l^{N_0}}\lambda_{k,l,l_1,l_2} e^{i\mu_l(x-a)},
\end{equation}
where the coefficients $\lambda_{k,l,l_1,l_2} $ are given by
\begin{align}
&\lambda_{k,l,l_1,l_2}=-\tau\mathcal{G}_{k,l,l_1,l_2}(\tau/2)+\int_0^\tau\mathcal{G}_{k,l,l_1,l_2}( s)\,d s
= r_{l,l_2}e^{it_k\delta_{l,l_2}}c_{k,l,l_1,l_2},
\label{eq:Lambdak1}
\end{align}
and
\begin{align}
&c_{k,l,l_1,l_2}= \widehat{V}_{l_1} \widehat{\phi}_{l_2}(t_k), \nn\\ 
&r_{l,l_2}=  \ -\tau e^{i\frac{\tau\delta_{l,l_2}}{2}}+\int_0^\tau e^{is\delta_{l,l_2}}\,ds = O(\tau^3 (\delta_{l,l_2})^2).
\label{eq:rest1}
\end{align}
The key observation from \eqref{eq:rest1} is: if $\delta_{l,l_2}=0$,  $r_{l,l_2}=0$ and the term $\lambda_{k,l,l_1,l_2}$ in \eqref{eq:Lambdak1} vanishes. Thus, in the discussion below, we shall assume that $\delta_{l,l_2}\neq0$.
Based on the RCO, we will go through the detailed structure of \eqref{eq:remainder-dec1} and exchange the order of summation (sum over index $k$ first), which will result in the terms like $\phi(t_k)-\phi(t_{k+1})=O(\tau\partial_t\phi)=O(\eps\tau)$ to gain an order of $\eps$.

First, for $l\in\mathcal{T}_{N_0}$ and $(l_1,l_2)\in\mathcal{I}_l^{N_0}$, we have
\begin{equation}
\vert \delta_{l,l_2}\vert \leq \delta_{N_0/2}= \mu_{{N_0}/2}^2=(\pi N_0)^2/(b-a)^2\leq\frac{4\pi^2 (1+\tau_0)^2}{\tau_0^2(b-a)^2}=\frac{\mu_1^2(1+\tau_0)^2}{\tau_0^2},
\end{equation}
which implies the following estimates for   the case \eqref{eq:tauc1} ($0< \tau \leq \alpha\frac{2\pi}{\mu_1^2(1+\tau_0)^2}\tau_0^2$ with $\alpha\in(0,1)$ and $\tau_0\in(0,1)$)
\begin{equation}
\label{eq:cfl1}
\frac{\tau}{2}\vert \delta_{l,l_2} \vert < \alpha\pi,
\end{equation}
and for the case \eqref{eq:tauc2}, we have 
\begin{equation}
\label{eq:cfl2}
\left|\frac{\tau}{2}\vert \delta_{l,l_2}\vert-k\pi \right| \ge\frac{\lambda}{2|\delta_{l,l_2}|^{1+\nu_3}},\quad k=\lceil\frac{|\delta_{l,l_2|\tau}}{2\pi}\rceil,\lceil\frac{|\delta_{l,l_2}|\tau}{2\pi}\rceil-1,
\end{equation}
where $\frac{\lambda}{2|\delta_{l,l_2}|^{1+\nu_3}}\leq \frac{\alpha\pi}{8}\leq\frac{\pi}{8}$ and $\sin\left(\frac{\lambda}{2|\delta_{l,l_2}|^{1+\nu_3}}\right)\ge \frac{\lambda}{4|\delta_{l,l_2}|^{1+\nu_3}}$ ($\sin(s)\ge\frac{s}{2}$ if $s\in(0,\pi/3)$).

Denoting $S_{n,l,l_2}=\sum\limits_{k=0}^n e^{it_k\delta_{l,l_2}}$ ($n\ge0$) and using summation by parts,
we find from \eqref{eq:Lambdak1} that
\begin{align}
\sum_{k=0}^n\lambda_{k,l,l_1,l_2} = & \
r_{l,l_2}\sum_{k=0}^{n-1}S_{k,l,l_2} (c_{k,l,l_1,l_2}-c_{k+1,l,l_1,l_2}) \nn \\
& \ +S_{n,l,l_2}\;r_{l,l_2}\;c_{n,l,l_1,l_2},\label{eq:lambdasum1}
\end{align}
and
\begin{align}
&c_{k,l,l_1,l_2}-c_{k+1,l,l_1,l_2} = \widehat{V}_{l_1}\left(\widehat{\phi}_{l_2}(t_k)-\widehat{\phi}_{l_2}(t_{k+1})\right).\label{eq:cksum1}
\end{align}
For the  step size \eqref{eq:tauc1}, we know from \eqref{eq:cfl1}  that
\begin{equation}
\label{eq:Sbd1}
 \vert S_{n,l,l_2}\vert \leq \frac{2}{\vert1-e^{i\tau\delta_{l,l_2}}\vert}=\frac{1}{\vert \sin(\tau \delta_{l,l_2}/2)\vert }\leq\frac{C}{\tau\vert\delta_{l,l_2}\vert},\quad\forall n\ge0,
\end{equation}
where  we have used the fact $\frac{\sin (s)}{s}$ is bounded (decreasing) for $s\in[0, \alpha\pi)$ and $C=\frac{2\alpha\pi}{\sin(\alpha\pi)}$ (noticing the case $\delta_{l,l_2}=0$ is trivial and $\delta_{l,l_2}$ is assumed to be nonzero here).
Combining \eqref{eq:rest1}, \eqref{eq:lambdasum1}, \eqref{eq:cksum1} and \eqref{eq:Sbd1}, we have
\begin{equation}
\left\vert\sum_{k=0}^n\lambda_{k,l,l_1,l_2}\right\vert
\lesssim   \tau^2\vert\delta_{l,l_2}\vert\left\vert\widehat{V}_{l_1}\right\vert\left[\sum\limits_{k=0}^{n-1} \left\vert\widehat{\phi}_{l_2}(t_k)-\widehat{\phi}_{l_2}(t_{k+1})\right\vert
 + \left\vert\widehat{\phi}_{l_2}(t_n)\right\vert\right] .\label{eq:sumlambda1}
\end{equation}
We note that \eqref{eq:sumlambda1} is the key for the refined estimate, where we shall gain an order of $\varepsilon$ from the $\phi(t_k)-\phi(t_{k+1})$ terms (see \eqref{eq:S1twist}). Of course, the condition \eqref{eq:cfl1} is also important to exclude the resonance case  where $S_n$ could be unbounded, i.e., \eqref{eq:cfl1} makes the estimate \eqref{eq:sumlambda1} available.  For the  non-resonance step size \eqref{eq:tauc2}, we can similarly obtain $ \vert S_{n,l,l_2}\vert \leq 1/\sin(\lambda/2|\delta_{l,l_2}|^{1+\nu_3})\leq \frac{4|\delta_{l,l_2}|^{1+\nu_3}}{\lambda}=\tilde{C}|\delta_{l,l_2}|^{1+\nu_3}/\alpha$ for some constant $\tilde{C}>0$.  Then, by noticing $\vert r_{l,l_2}\vert=O(\tau^2\delta_{l,l_2})$, the simialr estimates in \eqref{eq:sumlambda1}  hold  as
\begin{equation*}
\left\vert\sum_{k=0}^n\lambda_{k,l,l_1,l_2}\right\vert
\lesssim   \tau^2\vert\delta_{l,l_2}\vert^{2+\nu_3}\left\vert\widehat{V}_{l_1}\right\vert\left[\sum\limits_{k=0}^{n-1} \left\vert\widehat{\phi}_{l_2}(t_k)-\widehat{\phi}_{l_2}(t_{k+1})\right\vert
 + \left\vert\widehat{\phi}_{l_2}(t_n)\right\vert\right] .
\end{equation*}
 Since the rest arguments are almost the same for both step sizes \eqref{eq:tauc1} and \eqref{eq:tauc2} (regularities are different), we shall only treat the case \eqref{eq:tauc1} below. 

{\bf{Step 3.}} (Improved estimates) Now, we are ready to give the improved estimates. For $l\in\mathcal{T}_{N_0}$ and $(l_1,l_2)\in\mathcal{I}_l^{N_0}$, simple calculations show ($l=l_1+l_2$)
\begin{equation}
1+\mu_l^2\leq (1+\mu^2_{l_1})(1+\mu^2_{l_2}),\quad \vert \delta_{l,l_2}\vert \leq (1+\mu^2_{l_1})(1+\mu^2_{l_2}).
\label{eq:mlbd1}
\end{equation}
Based on \eqref{eq:remainder-dec1}, \eqref{eq:sumlambda1} and \eqref{eq:mlbd1},  using Cauchy inequality, we now estimate the remainder term in \eqref{eq:final1},
\begin{align}
&\left\| \mathcal{R}^n(x)\right\|^2_{H^1}  \label{eq:sumlambda-21}\\
&= \varepsilon^2
\sum\limits_{l\in\mathcal{T}_{N_0}} \left(1+\mu^2_l\right) \big\vert \sum\limits_{(l_1,l_2)\in\mathcal{I}_l^{N_0}}\sum\limits_{k=0}^n\lambda_{k,l,l_1,l_2}\big\vert^2\nonumber\\
&\lesssim \varepsilon^2 \tau^4
\bigg\{\sum_{l\in\mathcal{T}_{N_0}}\bigg(\sum\limits_{(l_1,l_2)\in\mathcal{I}_l^{N_0}}\left\vert \widehat{V}_{l_1}\right\vert \left\vert \widehat{\phi}_{l_2}(t_n)\right\vert \prod_{j=1}^2(1+\mu_{l_j}^2)^{3/2} \bigg)^{2}\nonumber\\
&\quad+ n \sum\limits_{k=0}^{n-1} \bigg[\sum_{l\in\mathcal{T}_{N_0}}\bigg(\sum\limits_{(l_1,l_2)\in\mathcal{I}_l^{N_0}}
\left\vert\widehat{V}_{l_1}\right\vert\left\vert\widehat{\phi}_{l_2}(t_k)-\widehat{\phi}_{l_2}(t_{k+1})\right\vert\prod_{j=1}^2(1+\mu_{l_j}^2)^{3/2} \bigg)^2\bigg]\bigg\}.\nonumber
\end{align}
To estimate each term in \eqref{eq:sumlambda-21}, we use the auxiliary function $\xi(x)=\sum_{l\in\mathbb{Z}}(1+\mu_l^2)^{3/2}\left\vert\widehat{\phi}_l(t_n)\right\vert e^{i\mu_l(x-a)}$, where $\xi(x)\in H_{\rm per}^{m-3}(\Omega)$ implied by the assumption (A) and $\|\xi(x)\|_{H^{s}}\lesssim\|\phi(t_n)\|_{H^{s+3}}$ ($s \leq m-3$). Similarly, introduce the function $U(x) = \sum_{l\in\mathbb{Z}}(1+\mu_l^2)^{3/2}\left\vert \widehat{V}_l\right\vert e^{i\mu_l(x-a)}$, where $U(x) \in H^{2}_{\rm per}$ implied by the assumption (B). Expanding
\[U(x)\xi(x)=\sum\limits_{l\in\mathbb{Z}}\sum\limits_{l_1+l_2=l} \prod_{j=1}^2 (1+\mu_{l_j}^2)^{3/2} \left\vert\widehat{V}_{l_1}\right\vert\left\vert\widehat{\phi}_{l_2}(t_n)\right\vert e^{i\mu_l(x-a)},\]
we could obtain
\begin{align}
&\sum_{l\in\mathcal{T}_{N_0}} \bigg(\sum\limits_{(l_1,l_2)\in\mathcal{I}_l^{N_0}}\left\vert\widehat{V}_{l_1}\right\vert\left\vert\widehat{\phi}_{l_2}(t_n)\right\vert \prod_{j=1}^2(1+\mu_{l_j}^2)^{3/2} \bigg)^2\nn\\
& \leq \ \|U(x)\xi(x)\|^2_{L^2} \lesssim \|V(x)\|_{H^4}^2 \|\phi(t_k)\|_{H^3}^2 \lesssim1,
\label{eq:est_ll}
\end{align}
which together with \eqref{eq:S1twist} implies (applying the same trick to the rest terms) for $0 \leq n \leq \frac{T/\varepsilon}{\tau}-1$,
\begin{align}
\left\| \mathcal{R}^n(x)\right\|^2_{H^1} & \lesssim \varepsilon^2\tau^4 \bigg(\|\phi(t_k)\|_{H^3}^2 + n\sum\limits_{k=0}^{n-1}
\|\phi(t_k) - \phi(t_{k+1})\|_{H^3}^2\bigg)\nn \\
& \lesssim \varepsilon^2 \tau^4+n^2\varepsilon^4 \tau^6 \|\partial_t \phi(x,t)\|_{L^\infty([0,T_\eps];H^3)}^2
\lesssim \varepsilon^2\tau^4.
\label{eq:est-l2l}
\end{align}
Combining \eqref{eq:final1} and \eqref{eq:est-l2l}, we have
\begin{equation}
\|e^{n+1}\|_{H^1} \lesssim h^{m-1}  + \tau_0^{m-1} + \varepsilon \tau^2 +\varepsilon\tau \sum\limits_{k=0}^n\|e^{k}\|_{H^1},\quad 0\leq n\leq\frac{T/\varepsilon}{\tau}-1.
\end{equation}
Discrete Gronwall's inequality would yield $\|e^{n+1}\|_{H^1} \lesssim h^{m-1} + \varepsilon \tau^2 + \tau_0^{m-1}$ $ (0\leq n\leq\frac{T/\varepsilon}{\tau}-1)$, and the proof for the improved uniform error bound \eqref{eq:error_l} in Theorem \ref{thm:im_linear} is completed.
\end{proof}

\begin{remark}\label{remark:hd} From the proof, the key steps of RCO  are the cut-off \eqref{eq:cut-off} to separate the high/low Fourier modes,   sufficiently small time step size $\tau$ \eqref{eq:tauc1} (or non-resonance step size \eqref{eq:tauc2}) to compensate the growth of errors at low Fourier modes via expansion (cf. \eqref{eq:cfl1},\eqref{eq:Sbd1} and \eqref{eq:sumlambda1}), and the estimates of the Fourier coefficients (cf. \eqref{eq:est-l2l}). Here, the special structure of the Fourier functions are important, e.g. $e^{i\mu_l(x-a)}e^{i\mu_k(x-a)}=e^{i\mu_{l+k}(x-a)}$. Based on above observations, it is straightforward to extend the RCO analysis to the higher dimensions (2D/3D) in rectangular domain with periodic boundary conditions for the sufficiently small time step sizes (\eqref{eq:tauc1} type), as the higher dimensional tensor Fourier basis enjoy the same properties ensuring the Fourier expansion for products of periodic functions.
 Notice that in 2D/3D for rectangular domains with irrational aspect ratios, the higher dimensional version of the non-resonance step size \eqref{eq:tauc2} is difficult to check, while the \eqref{eq:tauc1} type condition always holds for sufficiently small $\tau$.
\end{remark}
\begin{remark} In the proof, \eqref{eq:cfl1}, \eqref{eq:Sbd1} and \eqref{eq:sumlambda1} suggest the two order spatial regularity is regained from the summation by parts process indicated by $1/(\tau\delta_{l,l_2})$ ($\delta_{l,l_2}$ is roughly $\Delta$). Usually, such gained regularity will be lost when considering the other term of the summation by parts, i.e. the terms corresponding to $\phi(t_k)-\phi(t_{k+1})=O(\tau\partial_t\phi)$, where $\partial_t$ term will compensate the regularity gain from  $\delta_{l,l_2}$. However, as we have chosen a particularly designed twisted variable $\phi(t)$, there will be no regularity loss in $\phi(t_k)-\phi(t_{k+1})$, but with a gain of order $\eps$.
\end{remark}

\begin{remark}\label{rmk:tau0}
Passing $h\to0^+$ in Theorem \ref{thm:im_linear}, we can recover the estimates in the semi-discrete-in-time case, and it would be interesting to derive the estimates involving $\tau$ and $h$ only, i.e. $O(h^{m-1}+\varepsilon\tau^2)$ without the parameter $\tau_0$.  The following two cases are included:

1. Non-resonance $\tau$. Since the free Schr\"odiner operator $e^{it\Delta}$ is periodic in $t$,  we could impose the following Diophantine condition \cite{DG,SZ}: there exists $\gamma > 0$ and $\nu > 1$ such that
\begin{equation}\label{eq:resinf}
\left|\frac{1-e^{i\tau \mu_1^2K}}{\tau}\right|\geq \frac{\gamma}{|K|^\nu},	\quad \forall K\in\mathbb{Z},\; K\neq0,
\end{equation}
which is a common choice allowing $\tau\to0^+$. In particular, one can choose $\tau$ similar to \eqref{eq:tauc2} as \cite{SZ}
\begin{equation}\label{eq:setf}
\left\{\tau\in(0,1): \left|\tau-\frac{2l\pi}{\mu_1^2K}\right|\ge\frac{\lambda}{|l|^{\tilde{\nu}_1}|K|^{\tilde{\nu}_2}},  K,l\in\mathbb{Z},\; 1\leq l, |K|\right\},
\end{equation}
where $\tilde{\nu}_1\in[-1,1]$ and $\tilde{\nu}_1+\tilde{\nu}_2>2$, and $\lambda>0$ is a small constant.
The above set \eqref{eq:setf} is nowhere dense and hence not easy to verify in practice for a particular choice of $\tau$.
Once the non-resonance condition \eqref{eq:resinf} is satisfied, the improved estimates in Theorem \ref{thm:im_linear} hold and \eqref{eq:nonres} holds for any $\tau_0\in(0,1)$. Therefore, we can simply take the limit as $\tau_0\to0^+$ to derive the improved error bounds at $O(\varepsilon\tau^2+h^{m-1})$ for the above non-resonance step sizes.

2. Sufficiently small $\tau$ with the constraint $\tau\lesssim C_\varepsilon$ for some $\varepsilon$ dependent $C_\varepsilon$ (this type step size $\tau$ may not be included in \eqref{eq:resinf}). Following the proof of Theorem \eqref{thm:im_linear}, we can fixe $\tau_0=C\sqrt{\tau}$ for some constant $C$. By optimizing the error bounds $O(\varepsilon\tau^2+\tau_0^{m-1})=O(\varepsilon\tau^2+C^{m-1}\tau^{(m-1)/2})$, we find that  the $O(\eps\tau^2+h^{m-1})$ error bound would hold for sufficiently small $\tau\lesssim \varepsilon^{2/(m-5)}$ when $m>5$. If the exact solution $\psi(x,t)$ is sufficiently smooth with Fourier coefficient $\hat{\psi}_l(t)\sim O(e^{-c\vert l\vert})$ ($c>0$, $\vert l\vert\gg1$) decaying  exponentially fast, the projection error due to the $\tau_0$ cut-off would be $O(e^{-\tilde{C}/\sqrt{\tau}})$ and the error bounds become $O(\varepsilon\tau^2+h^{m-1})$ when $\tau\lesssim 1/\vert\ln\varepsilon\vert$ ($\varepsilon\in(0,1)$). Therefore, $C_\varepsilon=\varepsilon^{2/(m-5)}$ for sufficiently smooth solutions with  $m>5$ and $C_\varepsilon=1/\vert\ln\varepsilon\vert$ for the solutions with exponentially decaying Fourier coefficients.


\end{remark}

\subsection{Numerical results}
In this subsection, we present numerical results of the TSFP method for the long-time dynamics of the Schr\"odinger equation with $O(\eps)$-potential in 1D, up to the time $T_\eps=\frac{T}{\varepsilon}$.

First, we show an example to confirm that the uniform error bound in $L^2$-norm linearly grows with respect to  $T$. We choose the potential $V(x) = 5 \cos(2\pi x)$ and the $H^2_{\rm per}$ initial data as
\begin{equation}
\psi_0(x) = 5x^2(1-x)^2, \quad x \in [0, 1].
\end{equation}
The regularity is enough to ensure the uniform and the improved error bounds in $L^2$-norm.
The `exact' solution $\psi(x, t)$ is obtained numerically by the TSFP \eqref{eq:TSFP} with a very fine mesh size $h_e = 1/128$ and  time step size $\tau_e = 10^{-4}$. 
To quantify the error, we introduce the following error functions:
\begin{equation}
e_{L^2}(t_n) = \left\|\psi(x, t_n) - I_N \psi^n\right\|_{L^2}, \quad e_{H^1}(t_n) = \left\|\psi(x, t_n) - I_N \psi^n\right\|_{H^1},		
\end{equation}
and
\begin{equation*}
e_{L^2, \max}(t_n) = \max_{0 \leq q \leq n}e_{L^2}(t_q), \quad e_{H^1, \max}(t_n) = \max_{0 \leq q \leq n}e_{H^1}(t_q).	
\end{equation*}
In the rest of the paper, the spatial mesh size is always chosen sufficiently small and thus spatial errors can be ignored when considering the long time error growth and/or the temporal errors.
\begin{figure}[ht!]
\centerline{\includegraphics[width=12cm,height=5.5cm]{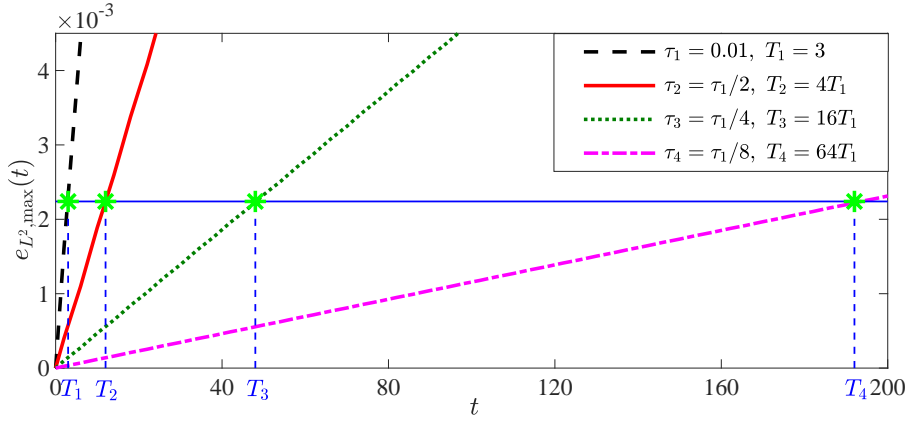}}
\caption{Long-time temporal errors in $L^2$-norm of the TSFP \eqref{eq:TSFP} for the Schr\"odinger equation \eqref{eq:linear_1D} with $\varepsilon = 1$ and different time step $\tau$.}
\label{fig:Long_linear_tau}
\end{figure}

\begin{figure}[ht!]
\centerline{\includegraphics[width=12cm,height=5.5cm]{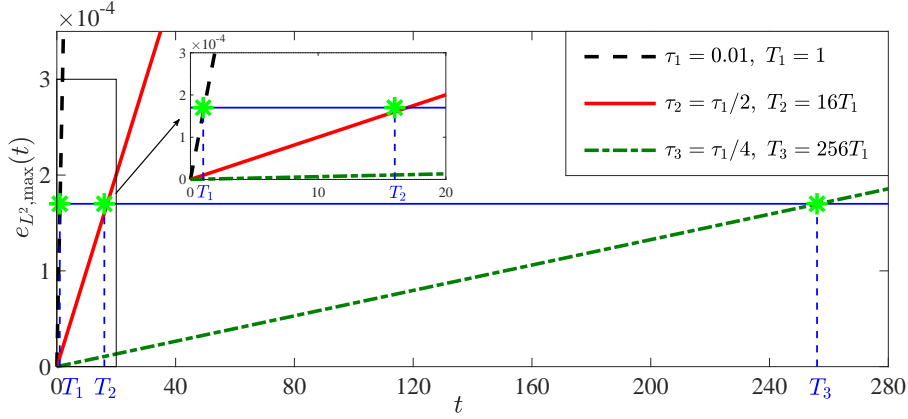}}
\caption{Long-time temporal errors in $L^2$-norm of the fourth-order time-splitting method for the Schr\"odinger equation \eqref{eq:linear_1D} with $\varepsilon = 1$ and different time step $\tau$.}
\label{fig:Long_linear_tau_4}
\end{figure}

Figure \ref{fig:Long_linear_tau} plots the long-time errors in $L^2$-norm of the TSFP method for the Schr\"odinger equation \eqref{eq:linear_1D} with $\varepsilon = 1$ and different time step $\tau$, which shows that the uniform errors in $L^2$-norm linearly grows with respect to the time. In addition, for a given accuracy bound, the time to exceed  the error bar is quadruple when the time step is half, which also confirms the linear growth. For comparisons, Figure \ref{fig:Long_linear_tau_4} depicts the long-time errors in $L^2$-norm of the fourth-order time-splitting method, which indicates that higher order time-splitting methods could get better accuracy with the same time step  size as well as longer time simulations within a given accuracy bound.

Next, we report the convergence test for the Schr\"odinger equation \eqref{eq:linear_1D} with the potential $V(x) = \sin(x)$ and the smooth initial data
\begin{equation}
\psi_0(x) = 2/(2+\sin^2(x)), \quad x \in [0, 2\pi].	
\label{eq:num_in}
\end{equation}
The `exact' solution $\psi(x, t)$ is obtained numerically by the TSFP \eqref{eq:TSFP} with  $h_e = \pi/64$ and  $\tau_e = 10^{-4}$.
\begin{figure}[ht!]
\centerline{\includegraphics[width=12cm,height=5.5cm]{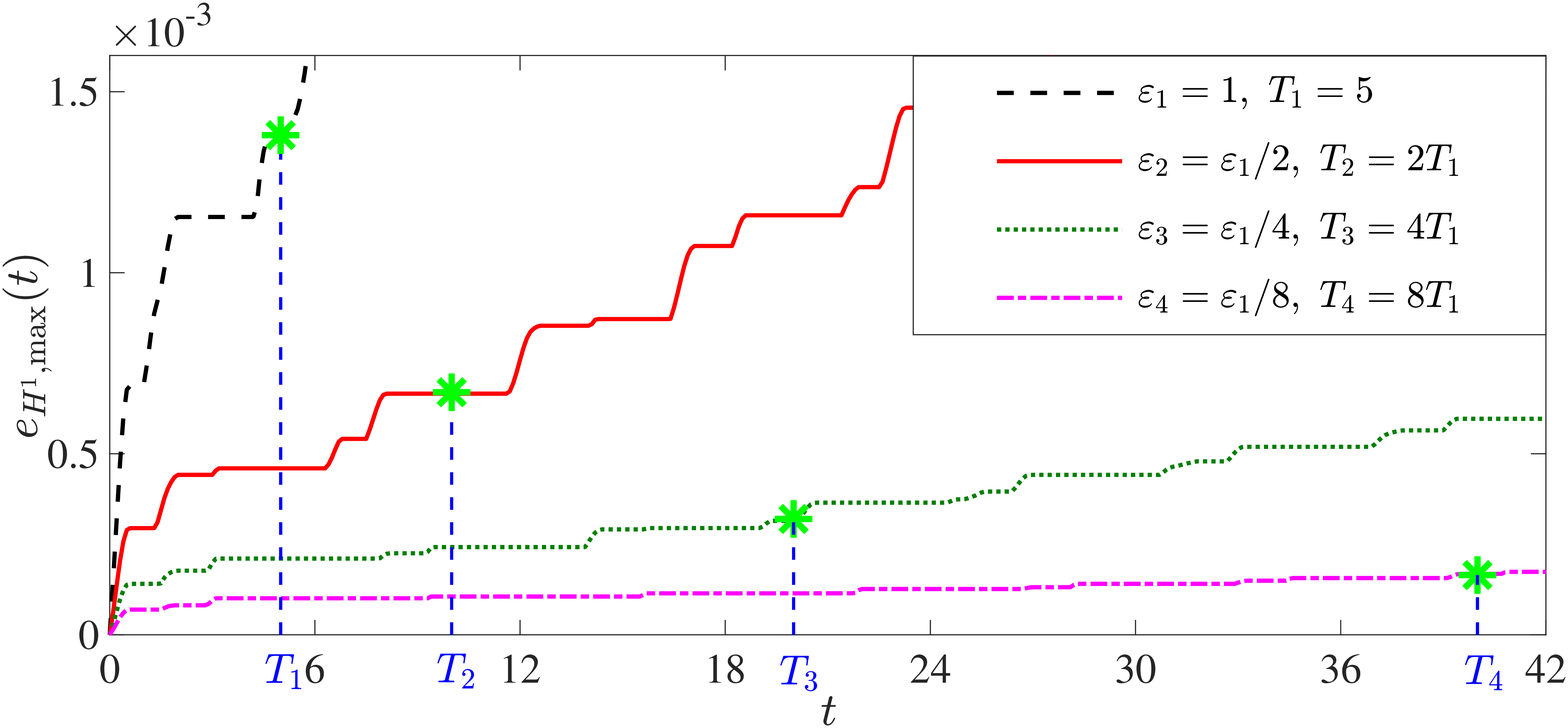}}
\caption{Long-time temporal errors in $H^1$-norm of the TSFP \eqref{eq:TSFP} for the Schr\"odinger equation \eqref{eq:linear_1D} with different $\varepsilon$.}
\label{fig:linear_long}
\end{figure}

\begin{figure}[ht!]
\begin{minipage}{0.49\textwidth}
\centerline{\includegraphics[width=6.5cm,height=5.5cm]{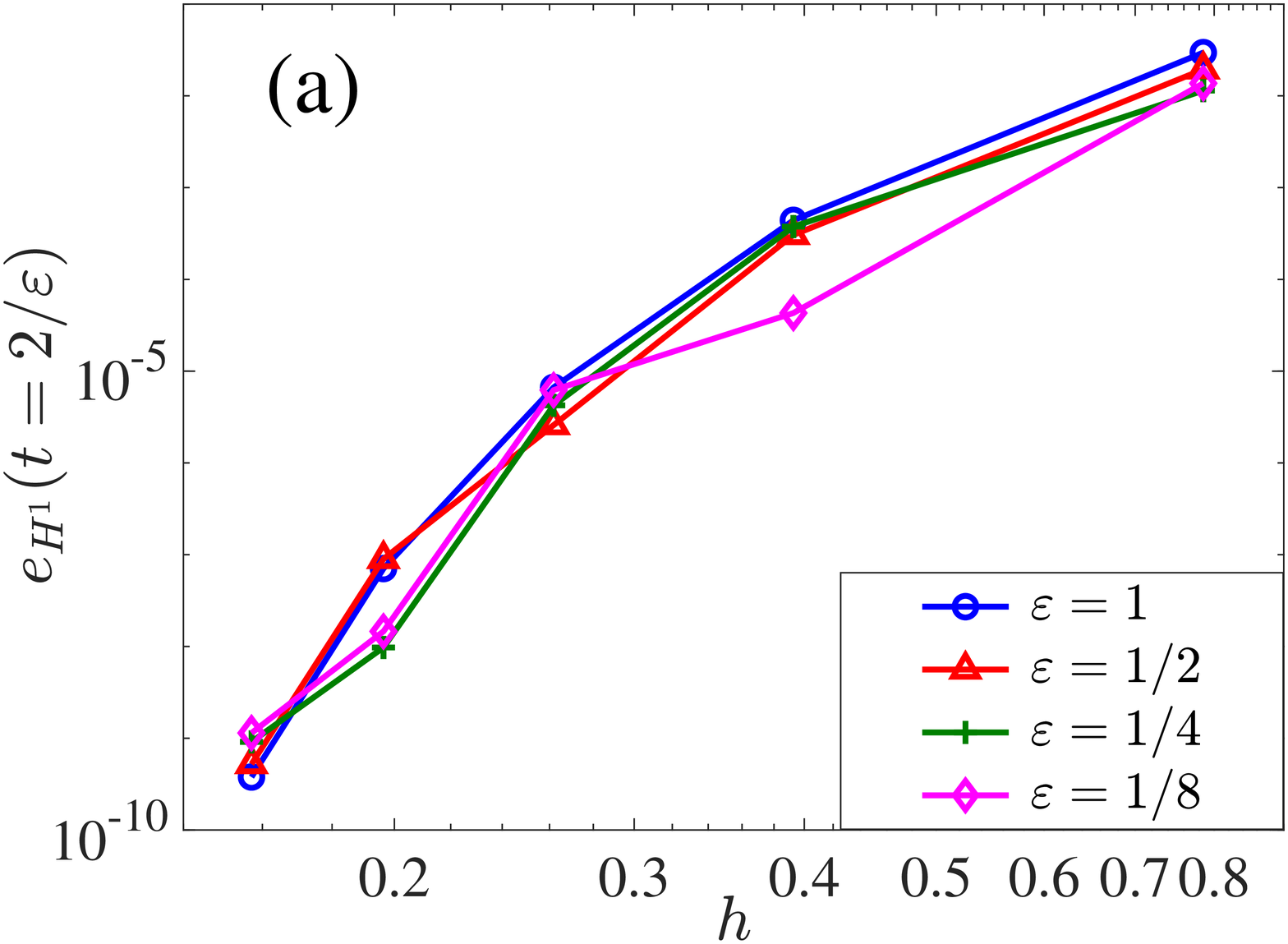}}
\end{minipage}
\begin{minipage}{0.49\textwidth}
\centerline{\includegraphics[width=6.5cm,height=5.5cm]{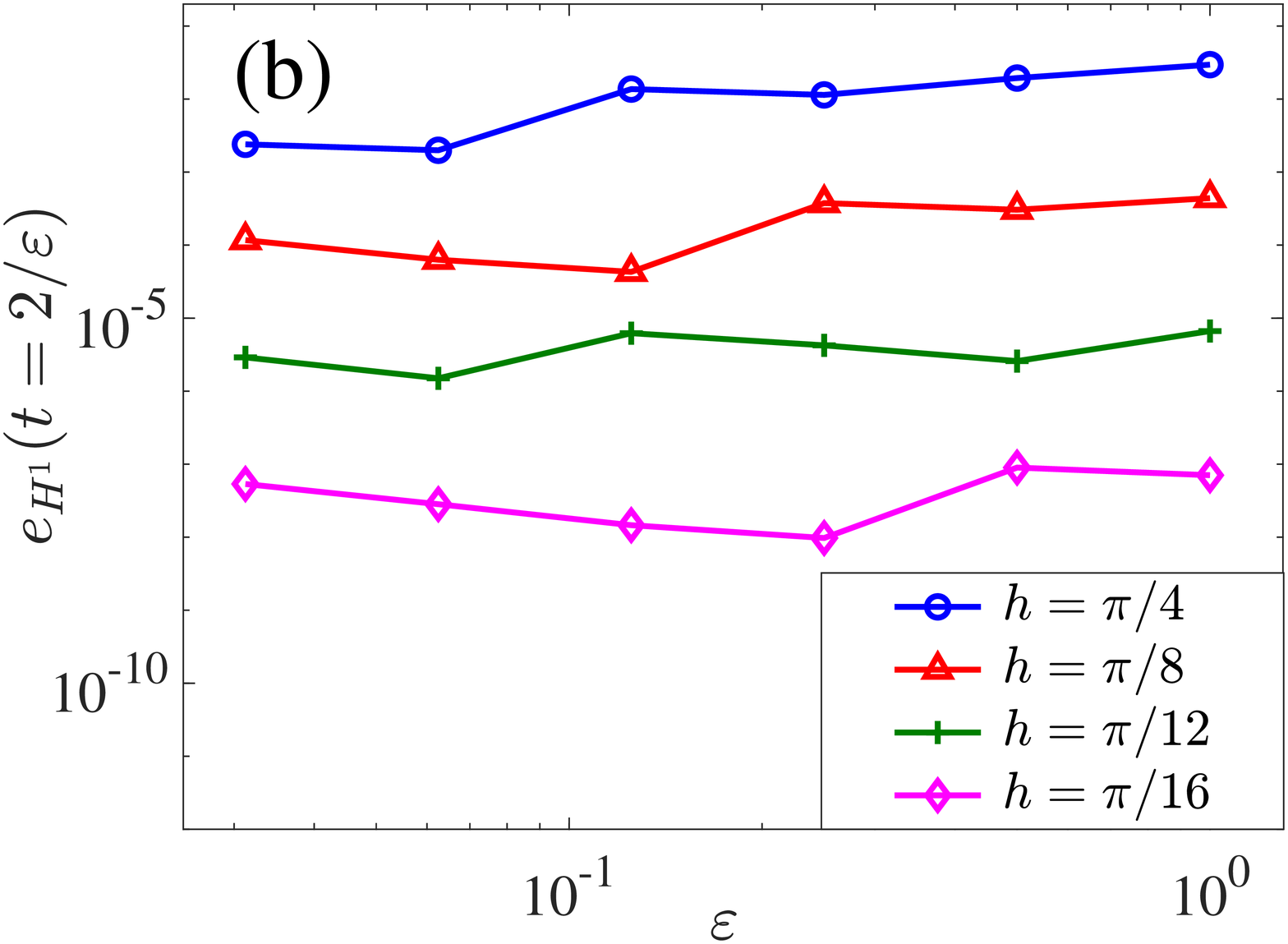}}
\end{minipage}
\caption{Long-time spatial errors in $H^1$-norm of the TSFP \eqref{eq:TSFP} for the Schr\"odinger equation \eqref{eq:linear_1D} at $t = 2/\varepsilon$.}
\label{fig:linear_spatial}
\end{figure}

\begin{figure}[ht!]
\begin{minipage}{0.49\textwidth}
\centerline{\includegraphics[width=6.3cm,height=5.5cm]{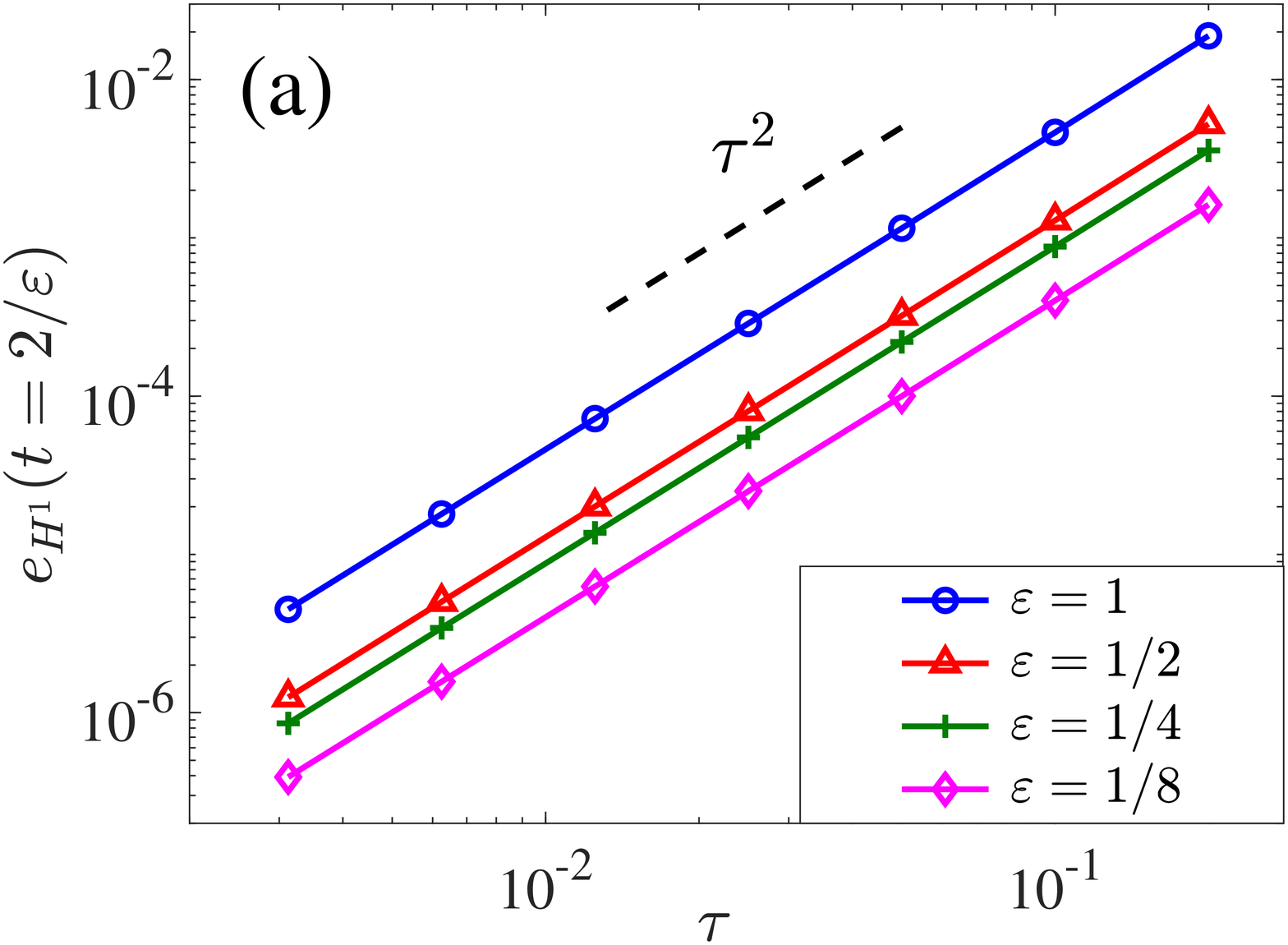}}
\end{minipage}
\begin{minipage}{0.49\textwidth}
\centerline{\includegraphics[width=6.3cm,height=5.5cm]{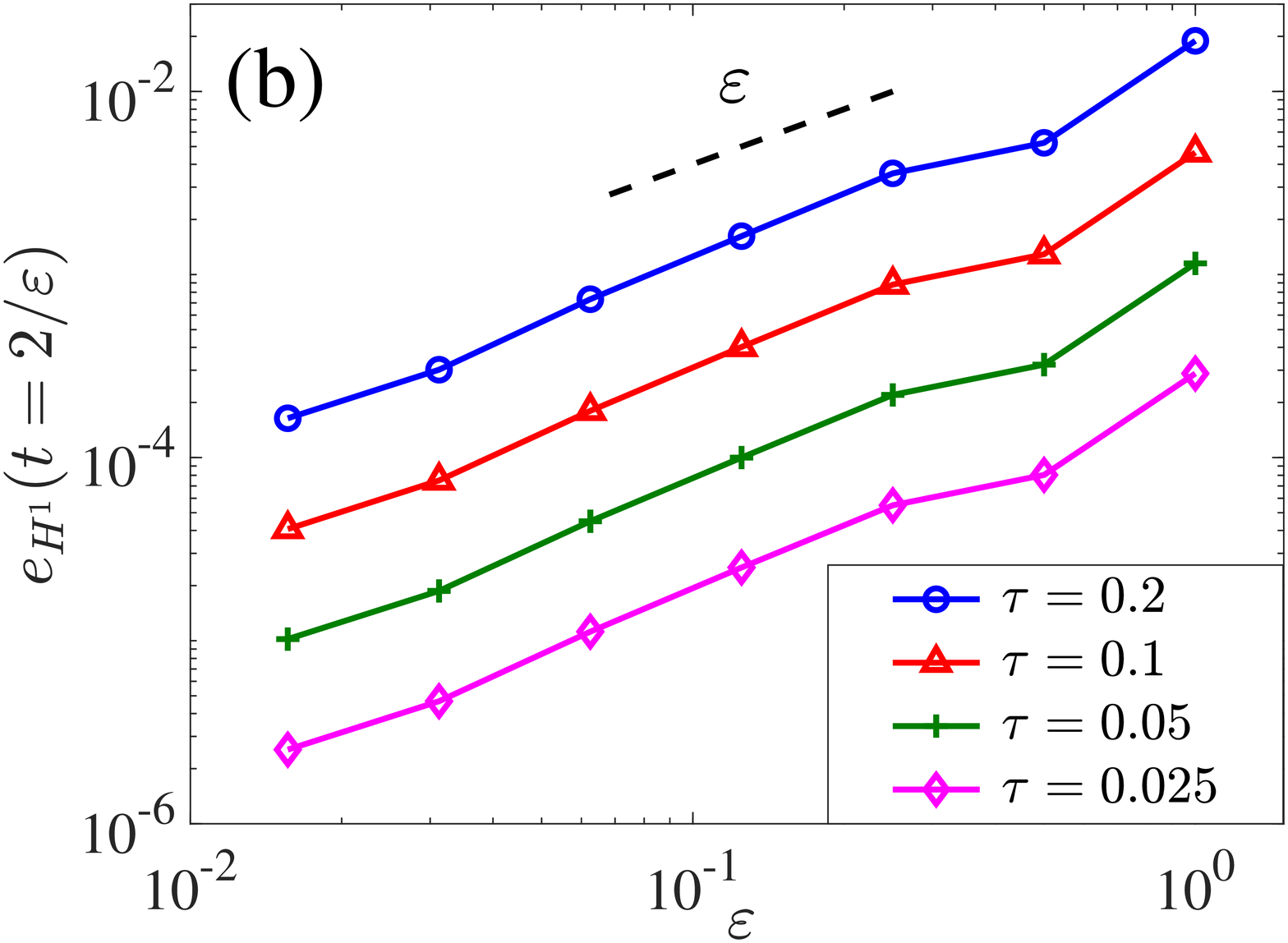}}
\end{minipage}
\caption{Long-time temporal errors in $H^1$-norm of the TSFP \eqref{eq:TSFP} for the Schr\"odinger equation \eqref{eq:linear_1D} at $t = 2/\varepsilon$.}
\label{fig:linear_temporal}
\end{figure}

Figure \ref{fig:linear_long} displays the long-time errors in $H^1$-norm of the TSFP method for the Schr\"odinger equation \eqref{eq:linear_1D} with the fixed time step $\tau$ and different $\varepsilon$, which confirms the improved uniform error bound in $H^1$-norm at $O(\eps \tau^2$) up to the $O(1/\eps)$ time. Figs. \ref{fig:linear_spatial} \&  \ref{fig:linear_temporal} exhibit the spatial and temporal errors of the TSFP \eqref{eq:TSFP} for the Schr\"odinger equation \eqref{eq:linear_1D}  at $t= \frac{2}{\varepsilon}$. Each line in Figure \ref{fig:linear_spatial} (a) shows the spectral accuracy of the TSFP method in space and Figure \ref{fig:linear_spatial} (b) verifies the spatial errors are independent of the small parameter $\varepsilon$ in the long-time regime. Figure \ref{fig:linear_temporal} (a) shows the second-order convergence of the TSFP method in time. Each line in Figure \ref{fig:linear_temporal} (b) gives the global errors in $H^1$-norm with a fixed time step $\tau$ and verifies that the global error performs like $O(\varepsilon\tau^2)$ up to the $O(1/\varepsilon)$ time.

\begin{figure}[ht!]
\begin{minipage}{0.49\textwidth}
\centerline{\includegraphics[width=6.3cm,height=5.5cm]{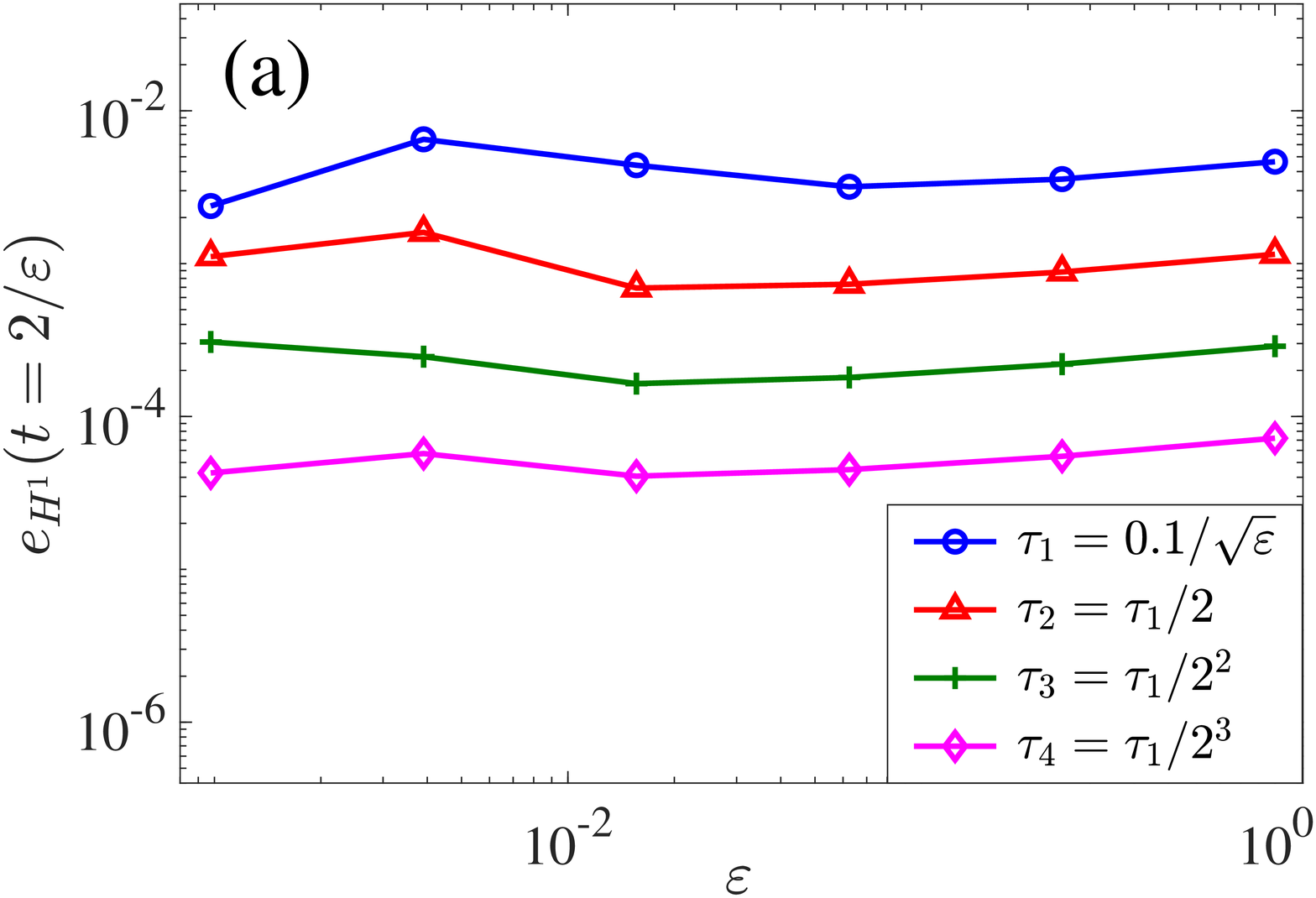}}
\end{minipage}
\begin{minipage}{0.49\textwidth}
\centerline{\includegraphics[width=6.3cm,height=5.5cm]{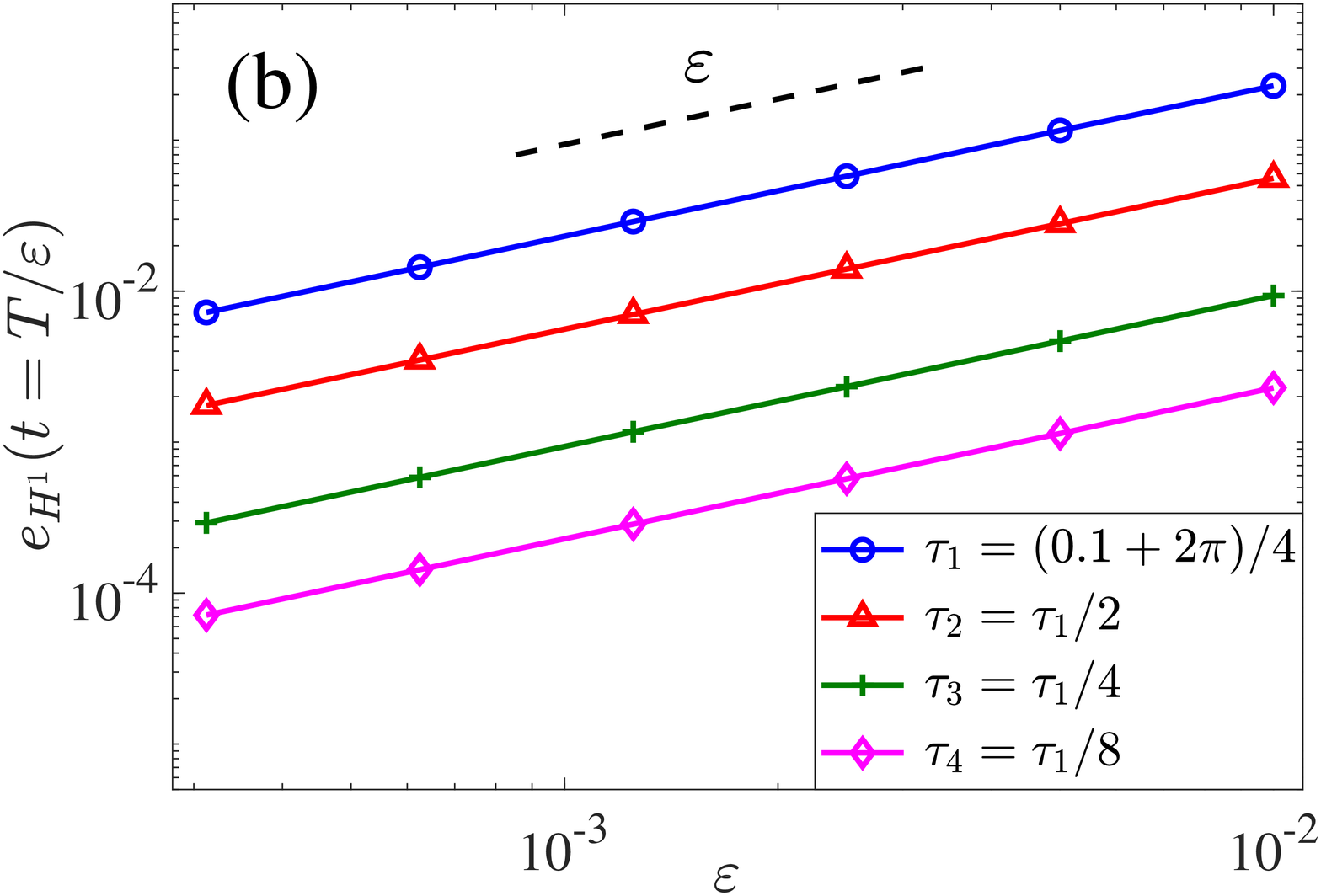}}
\end{minipage}
\caption{Long-time temporal errors in $H^1$-norm of the TSFP \eqref{eq:TSFP} for the Schr\"odinger equation \eqref{eq:linear_1D} with large time step size.}
\label{fig:linear2}
\end{figure}

Figure \ref{fig:linear2} displays the long-time errors of the TSFP method for the Schr\"odinger equation \eqref{eq:linear_1D} with large time step size. Each lines in Figure \ref{fig:linear2} (a) plots the long-time errors with the time step sizes $\tau = O(1/\sqrt{\eps})$, which are almost  constants for different $\eps$, confirming the error bound \eqref{eq:error_2}. In Figure \ref{fig:linear2} (b), we choose the time step size  satisfying the non-resonance condition, and Figure \ref{fig:linear2} (b) shows the improved uniform error bounds with large time step size.

\section{Improved uniform error bounds for the NLSE}
In this section, we adopt the TSFP method to solve the NLSE with weak nonlinearity and extend the technique of regularity compensation oscillation (RCO) to obtain improved uniform error bounds for the cubic NLSE with $O(\eps^2)$-nonlinearity up to the $O(1/\eps^2)$ time.

\subsection{The TSFP method}
We  present the TSFP method for the NLSE \eqref{eq:WNE} in 1D and extensions to higher dimensions are straightforward (see also Remark \ref{remark:hd}). In 1D, the NLSE \eqref{eq:WNE} with initial data \eqref{eq:initial} and periodic boundary conditions on $\Omega = (a, b)$ collapses to
\begin{equation}
\left\{
\begin{split}
& i\partial_t \psi(x, t) = -\Delta \psi(x, t) \pm \eps^2\vert\psi(x, t)\vert^2\psi(x, t),\quad a<x<b,\ t > 0, \\
&\psi(a, t) = \psi(b, t),\ \partial_x \psi(a, t) = \partial_x \psi(b, t),\ t \geq 0,\\
&\psi(x, 0) = \psi_0(x), \ x \in [a, b].
\end{split}\right.
\label{eq:nl_1D}
\end{equation}

By  the same time-splitting technique as that in the linear case, the semi-discretization of the NLSE \eqref{eq:nl_1D} via the Strang splitting is given as:
\begin{equation}
\psi^{[n+1]}(x) = \mathcal{S}_{\tau}(\psi^{[n]})= e^{i\frac{\tau}{2}\Delta}e^{\mp i \eps^2 \tau  \left\vert e^{i\frac{\tau}{2}\Delta}\psi^{[n]}(x)\right\vert^2} e^{i\frac{\tau}{2}\Delta}\psi^{[n]}(x), \ x \in \Omega,
\end{equation}
with $ \psi^{[0]}(x) = \psi_0(x)$. Respectively, the full-discretization for the NLSE \eqref{eq:nl_1D} can be written as
\begin{equation}
\label{eq:TSFP_nl}
\begin{split}
&\psi^{(1)}_j=\sum_{l \in \mathcal{T}_N} e^{-i\frac{\tau\mu^2_l}{2}}\;\widetilde{(\psi^n)}_l\; e^{i\mu_l(x_j-a)}, \\
&\psi^{(2)}_j=e^{\mp i\eps^2 \tau \lambda \vert\psi^{(1)}_j\vert^2} \psi^{(1)}_j, \qquad \qquad  j\in \mathcal{T}^0_N, \quad n \ge 0, \\
&\psi^{n+1}_j=\sum_{l \in \mathcal{T}_N} e^{-i\frac{\tau\mu^2_l}{2}} \; \widetilde{\left(\psi^{(2)}\right)}_l\; e^{i\mu_l(x_j-a)},
\end{split}
\end{equation}
where $\psi^0_j = \psi_0(x_j)$ for $j\in \mathcal{T}^0_N$.

\subsection{Improved uniform error bounds in $H^1$-norm}
For the NLSE, we assume the exact solution $\psi(x, t)$ up to the time at $T_{\eps} = T/\eps^2$ with $T > 0$ fixed satisfies:
\[
{\rm(C)} \qquad
 \left\|\psi(x, t)\right\|_{L^{\infty}\left([0, T_{\varepsilon}]; H^{m}_{\rm per}\right)} \lesssim 1, \quad  \left\|\partial_t \psi(x, t) \right\|_{L^{\infty}\left([0, T_{\varepsilon}]; H^{m-2}_{\rm per}\right)} \lesssim 1,\quad m\ge5.
\]
Similar to Theorem \ref{thm:im_linear} in the linear case, we shall impose the following non-resonance conditions on the step size $\tau$ for TSFP \eqref{eq:TSFP_nl}  in the nonlinear case. For the Fourier modes $|l|\leq\lceil\frac{1}{\tau_0}\rceil$ ($\tau_0\in(0,1)$), we impose the condition: there exists a constant  $C_0>0$ such that
\begin{equation}\label{eq:nonres-nl}
\left|1-e^{i\tau \mu_1^2 K}\right|\ge \frac{C_0\tau^{\nu_1}}{(\mu_1^2 |K|)^{\nu_2}},\quad 0<|K|\leq K_1=2\lceil1/\tau_0\rceil^2,\quad K\in\mathbb{Z},
\end{equation}
where $\nu_1\in[0,1]$, $\nu_2\ge -1$, and the bound $|K|\leq 2\lceil1/\tau_0\rceil^2$ corresponds to the cubic nonlinear interaction. In particular, we consider the following cases of time step sizes: for a given constant $\alpha\in(0,1)$, the time step size $\tau$ satisfies
\begin{equation}\label{eq:tauc1-nl}
\tau\in\left(0, \alpha\frac{\pi}{\mu_1^2(1+\tau_0)^2}\tau_0^2\right),\end{equation}
or
\begin{equation}\label{eq:tauc2-nl}
\tau\in \left\{\tau>0: \left|\tau-\frac{2l\pi}{\mu_1^2K}\right|\ge\frac{\lambda}{|\mu_1^2K|^{2+\nu_3}},  K,l\in\mathbb{Z},\; 0<|K|\leq K_1,\;
0\leq l\right\},
\end{equation}
where $\lambda$ and $\nu_3$ are the same as those in \eqref{eq:tauc2}.

 Then we have the following improved uniform error bound of the TSFP \eqref{eq:TSFP_nl} for the NLSE with $O(\varepsilon^2)$-nonlinearity strength up to the time at $O(1/\varepsilon^2)$.

\begin{theorem}
\label{thm:eb_nl}
Let $\psi^n$ be the numerical approximation obtained from the TSFP \eqref{eq:TSFP_nl}. Under the assumption (C), there exist $h_0 > 0$, $0 < \tau_0 < 1$ sufficiently small and independent of $\varepsilon$ such that, for any $0 < \varepsilon \leq 1$, when $0 < h < h_0$ and 
 $\tau$ satisfies \eqref{eq:tauc1-nl} or \eqref{eq:tauc2-nl} ($m\ge5+\nu_3$ for \eqref{eq:tauc2-nl}) with $0<\tau\leq \tau_1/\varepsilon$ ($\tau_1>0$ small enough independent of $\varepsilon$)
the following error bounds hold
\begin{equation}
\begin{split}
&\left\|\psi(x, t_n) - I_N \psi^n\right\|_{H^1} \lesssim h^{m-1} + \varepsilon^2\tau^2 + \tau_0^{m-1}, \\
&\left\|I_N\psi^n\right\|_{H^1} \leq 1 + M, \ 0 \leq n \leq \frac{T/\varepsilon^2}{\tau},
\end{split}
\label{eq:error_nl}
\end{equation}
where $M := \|\psi\|_{L^{\infty}([0, T_{\varepsilon}]; H^1)}$. In particular, if the exact solution is smooth, i.e.  $\psi(x, t) \in H^{\infty}_{\rm per}$,  the $\tau_0^{m-1}$ error part would decrease exponentially and can be ignored in practical computation when $\tau_0$ is small but fixed, and thus the estimate would practically become
\begin{equation}
\left\|\psi(x, t_n) - I_N \psi^n\right\|_{H^1} \lesssim h^{m-1} + \varepsilon^2\tau^2.
\end{equation}
\end{theorem}
\begin{remark}Analogous to the linear case, improved error bounds \eqref{eq:error_nl} in Theorem \ref{thm:eb_nl} can be easily generalized to the non-resonance step size \eqref{eq:nonres-nl} under slightly different regularity assumptions. We also need  $\tau\leq \tau_1/\varepsilon$ for controlling the nonlinearity.
\end{remark}

Some parts of the proof proceed in analogous lines as the linear case and we omit the details in this section for brevity. Similar to the analysis of the local truncation error for the linear case, we have the following results for the local truncation error for the TSFP \eqref{eq:TSFP_nl}.
\begin{lemma}
 The local truncation error of the TSFP method \eqref{eq:TSFP_nl} for the NLSE with $O(\eps^2)$-nonlinearity strength can be written as $(0 \leq n \leq \frac{T/\eps^2}{\tau}-1)$
\begin{equation}\label{eq:local-nonlinear}
\overline{\mathcal{E}}^{n} :=P_N \mathcal{S}_{\tau}(P_N\psi(t_n)) - P_N\psi(t_{n+1}) = P_N\mathcal{J}(P_N\psi(t_n)) + Y_{n},
\end{equation}
where
\begin{equation}
\mathcal{J}(P_N\psi(t_n)) = -i\eps^2\tau g\left(\frac{\tau}{2}\right) + i\eps^2\int^{\tau}_0 g(s) ds,	
\label{eq:lambda_nl}
\end{equation}
with
\begin{equation}
 g(s) = \pm e^{i(\tau-s)\Delta} \vert P_N(\psi(t_n + s))\vert^2 e^{is\Delta}P_N\psi(t_n).
 \label{eq:gdef}
\end{equation}
Under the assumption (D), for $0 < \eps \leq 1$, we have the error bounds
\begin{equation}
\left\|\mathcal{J}(P_N\psi(t_n))\right\|_{H^1}\lesssim \eps^2 \tau^3\|\psi(t_n)\|_{H^5}^3,\quad \left\|Y_{n}\right\|_{H^1} \lesssim \varepsilon^4\tau^3+\eps^2\tau h^{m-1}.
\label{eq:local_nl}
\end{equation}
\label{lemma_nl}
\end{lemma}

\noindent
\emph{Proof for Theorem \ref{thm:eb_nl}.} We apply a standard induction argument for proving \eqref{eq:error_nl}. Since $\psi^0_j = \psi_0(x_j)$, it is obvious for $n = 0$. Assuming the error bounds \eqref{eq:error_nl} hold true for all $0 \leq n \leq q \leq \frac{T/\varepsilon^2}{\tau}-1$, we are going to prove the case $n = q + 1$. By Fourier projections $\|\psi(x, t_n)-I_N\psi^n\|_{H^1}\lesssim \|P_N\psi(x,t_n)-I_N\psi^n\|_{H^1}+h^{m-1}$, we  just need to analyze the growth of the error $e^n=I^n\psi^n - P_N\psi(t_n)$  carefully. For $0 \leq n \leq q$, we have
\begin{align}
e^{n+1} =  & \ I_N\psi^{n+1} - \mathcal{S}_{\tau}(P_N\psi(t_n)) + \overline{\mathcal{E}}^{n}= e^{i\tau\Delta}e^n +Z^n(x) + \overline{\mathcal{E}}^{n},
 \label{eq:etg_nl}
\end{align}
where $Z^n(x)$ is given by
\begin{align*}
Z^n(x) = &e^{i\frac{\tau}{2}\Delta}\bigg[I_N((e^{-i\eps^2\tau \lambda \vert e^{i\frac{\tau}{2}\Delta}I_N\psi^n\vert^2}-1)e^{i\frac{\tau}{2}\Delta}I_N\psi^{n})\\
& \quad\quad\quad - P_N((e^{-i\eps^2\tau \lambda\vert e^{i\frac{\tau}{2}\Delta} P_N\psi(t_n)\vert^2}-1)e^{i\frac{\tau}{2}\Delta}P_N\psi(t_n))\bigg],
\end{align*}
with the bound (constant in front of $\|e^n\|_{H^1}$ depends on $M$ )
\begin{equation}
\left\|Z^n(x)\right\|_{H^1} \lesssim \eps^2 \tau \left(h^{m-1} + \|e^n\|_{H^1}\right).
\end{equation}
From \eqref{eq:etg_nl},  we obtain for $0\leq n\leq q$,
\begin{align}\label{eq:n}
e^{n+1}=e^{i (n+1) \tau\Delta}e^0+\sum\limits_{k = 0}^n e^{i(n-k)\tau \Delta}\left(Z^k(x)+ {\overline{\mathcal{E}}}^{k}\right).
\end{align}
Similar to the linear case,  we get for $0 \leq n \leq q$,
\begin{align}
\|e^{n+1}\|_{H^1} \lesssim & \ h^{m-1}+\varepsilon^2\tau^2+\varepsilon^2\tau\sum_{k=0}^n\|e^k\|_{H^1}\nn \\
&+\left\| \sum\limits_{k=0}^n e^{i(n-k)\tau\Delta}P_N\mathcal{J}(P_N\psi(t_k))\right\|_{H^1}.
\label{eq:final}
\end{align}
Recalling \eqref{eq:lambda_nl} and \eqref{eq:gdef}, we could decompose $\mathcal{J}(\psi(t_n))$ as
\begin{equation}
\mathcal{J}(P_N\psi(t_n)) = \mathcal{J}_1(P_N\psi(t_n)) + \mathcal{J}_2(P_N\psi(t_n)),
\end{equation}
where $\mathcal{J}_{\sigma}(P_N\psi(t_n)) = -i\eps^2\tau g_{\sigma}(\tau/2) + i\eps^2\int^{\tau}_0 g_{\sigma}(s) ds$ for $\sigma = 1, 2$ and $g_{\sigma}(s) := g_{\sigma}(s; P_N\psi(t_n))$ $ (\sigma = 1, 2)$ are defined as
\begin{equation*}
 g_1(s) = \pm e^{i(\tau-s)\Delta} \vert e^{is\Delta}P_N\psi(t_n)\vert^2 e^{is\Delta}P_N\psi(t_n), \quad g_2(s) = g(s) - g_1(s),	
\end{equation*}
 with $g(s)$ $(s \in [0, \tau])$ given in \eqref{eq:gdef}. Under the assumption (D), by the Duhamel's principle, it is easy to verify $\|\vert e^{is\Delta}P_N\psi(t_n)\vert^2 - \vert P_N\psi(t_n+s)\vert^2\|_{L^{\infty}([0, \tau]; H^m)} \lesssim \eps^2 \tau$. Following similar analysis for the local truncation error in Section 2,  for $ 0 \leq n \leq \frac{T/\eps^2}{\tau}-1$, we could arrive at
 \begin{equation}
 \left\|\mathcal{J}_1(P_N\psi(t_n))\right\|_{H^1} \lesssim \eps^2\tau^3,	\quad  \left\|\mathcal{J}_2(P_N\psi(t_n))\right\|_{H^1} \lesssim \eps^4\tau^3.
 \end{equation}
In light of \eqref{eq:final}, we find the major part of the error is from $\mathcal{J}_1(\psi_n)$.

Following the RCO approach in the linear case, we introduce the  `twisted variable' $\phi(x, t)=e^{-it\Delta}\psi(x,t)$, and $\|\partial_t\phi\|_{L^\infty([0,T/\eps^2];H^m)}\lesssim \eps^2$ with
\begin{equation}
\left\|\phi(t_n) - \phi(t_{n-1})\right\|_{H^{m}} \lesssim \eps^2 \tau, \quad 1 \leq n \leq \frac{T/\eps^2}{\tau}.
\label{eq:nl_twist}
\end{equation}
We choose the same cut-off parameter $\tau_0\in(0,1)$ and the corresponding Fourier modes $N_0=2\lceil 1/\tau_0\rceil$ as in the proof of Theorem \ref{thm:im_linear}.
Thus, we can derive
\begin{align}
\left\|e^{n+1}\right\|_{H^1} \lesssim & \ h^{m-1} + \tau_0^{m-1} +\eps^2\tau^2+\eps^2\tau\sum_{k=0}^n\left\|e^k\right\|_{H^1} +\|\mathcal{L}^n\|_{H^1}, \label{eq:final2}\\
\mathcal{L}^n(x)=&\sum\limits_{k=0}^n e^{-i(k+1)\tau\Delta}P_{N_0}\mathcal{J}_1(e^{i t_k\Delta}(P_{N_0}\phi(t_k))).
\end{align}
For $l\in\mathcal{T}_{N_0}$, we define the index set $\mathcal{I}_l^{N_0}$ associated to $l$ as
\begin{equation}
\mathcal{I}_l^{N_0}=\left\{(l_1,l_2,l_3)\ \vert \ l_1-l_2+l_3=l,\ l_1,l_2,l_3\in\mathcal{T}_{N_0}\right\}.
\end{equation}
Then, the expansion below follows
\begin{equation*}
e^{-i(k+1)\tau\Delta}P_{N_0}\mathcal{J}_1(e^{i t_k\Delta}(P_{N_0}\phi(t_k))) =\sum\limits_{l\in\mathcal{T}_{N_0}}\sum\limits_{(l_1,l_2,l_3)\in\mathcal{I}_l^{N_0}}
\mathcal{G}_{k,l,l_1,l_2,l_3}(s)e^{i\mu_l(x-a)},
\end{equation*}
where the coefficients $\mathcal{G}_{k,l,l_1,l_2,l_3}(s)$ are functions of $s$ only,
\begin{equation}
\mathcal{G}_{k,l,l_1,l_2,l_3}(s) = e^{i(t_k+s)\delta_{l,l_1,l_2,l_3}}\left(\widehat{\phi}_{l_2}(t_k)\right)^{\ast}\widehat{\phi}_{l_1}(t_k)\widehat{\phi}_{l_3}(t_k),
\label{eq:mGdeff}
\end{equation}
and $\delta_{l,l_1,l_2,l_3} = \delta_l - \delta_{l_1} + \delta_{l_2} - \delta_{l_3}$ ($\delta_l=\mu_l^2$ as in \eqref{eq:mGdeff}).
The remainder term in \eqref{eq:final} reads
\begin{align}\label{eq:remainder-dec}
\mathcal{L}^n(x)  = \pm i\varepsilon^2\sum\limits_{k=0}^n
\sum\limits_{l\in\mathcal{T}_{N_0}}\sum\limits_{(l_1,l_2,l_3)\in\mathcal{I}_l^{N_0}}\Lambda_{k,l,l_1,l_2,l_3} e^{i\mu_l(x-a)},
\end{align}
where
\begin{align}
\Lambda_{k,l,l_1,l_2,l_3} &= - \tau \mathcal{G}_{k,l,l_1,l_2,l_3}(\tau/2) +\int_0^\tau\mathcal{G}_{k,l,l_1,l_2,l_3}( s)\,d s\nn\\
 &= r_{l,l_1,l_2,l_3}e^{it_k\delta_{l,l_1,l_2,l_3}}c_{k,l,l_1,l_2,l_3},
\label{eq:Lambdak}
\end{align}
with coefficients  $c_{k,l,l_1,l_2,l_3}$ and $r_{l,l_1,l_2,l_3}$ given by
\begin{align}
c_{k,l,l_1,l_2,l_3}=& \ (\widehat{\phi}_{l_2}(t_k))^{\ast}\widehat{\phi}_{l_1}(t_k)\widehat{\phi}_{l_3}(t_k),\label{eq:calFdef}\\
r_{l,l_1,l_2,l_3}= & \ -\tau e^{i\tau\delta_{l,l_1,l_2,l_3}/2}+\int_0^\tau e^{is\delta_{l,l_1,l_2,l_3}}\,d s \nn \\
= & \ O\left(\tau^3 (\delta_{l,l_1,l_2,l_3})^2\right).\label{eq:rest}
\end{align}
Similar to the linear case, we only need consider the case $\delta_{l,l_1,l_2,l_3}\neq0$,  as $r_{l,l_1,l_2,l_3}=0$ if $\delta_{l,l_1,l_2,l_3}=0$ . First, for $l\in\mathcal{T}_{N_0}$ and $(l_1,l_2,l_3)\in\mathcal{I}_l^{N_0}$, we have
\begin{equation}
\vert\delta_{l,l_1,l_2,l_3}\vert \leq 2\delta_{N_0/2}= 2 \mu_{N_0/2}^2 \leq \frac{8\pi^2 (1+\tau_0)^2}{\tau_0^2(b-a)^2}=\frac{2(1+\tau_0)^2\mu_1^2}{\tau_0^2},
\end{equation}
which implies for the case \eqref{eq:tauc1-nl} $0< \tau \leq \alpha \frac{\pi \tau_0^2}{\mu_1^2 (1+\tau_0)^2}$ with $0<\tau_0, \alpha<1$,
\begin{equation}\label{eq:cfl}
\frac{\tau}{2}\vert \delta_{l,l_1,l_2,l_3} \vert \leq \alpha\pi.
\end{equation}
Denoting $S_{n,l,l_1,l_2,l_3}=\sum_{k=0}^ne^{it_k\delta_{l,l_1,l_2,l_3}}$ ($n\ge0$) and using summation by parts,
we find from \eqref{eq:Lambdak} that
\begin{align}
\sum_{k=0}^n\Lambda_{k,l,l_1,l_2,l_3} =& \
r_{l,l_1,l_2,l_3}\sum_{k=0}^{n-1}S_{k,l,l_1,l_2,l_3} \left(c_{k,l,l_1,l_2,l_3}-c_{k+1,l,l_1,l_2,l_3}\right)\nonumber\\
&\ +S_{n,l,l_1,l_2,l_3}\, r_{l,l_1,l_2,l_3}\,c_{n,l,l_1,l_2,l_3},\label{eq:lambdasum}
\end{align}
and
\begin{align}
&c_{k,l,l_1,l_2,l_3}-c_{k+1,l,l_1,l_2,l_3}\nn\\
& = (\widehat{\phi}_{l_2}(t_k))^{\ast}(\widehat{\phi}_{l_1}(t_k)-\widehat{\phi}_{l_1}(t_{k+1})) \widehat{\phi}_{l_3}(t_k) + (\widehat{\phi}_{l_2}(t_k)-\widehat{\phi}_{l_2}(t_{k+1}))^{\ast}\widehat{\phi}_{l_1}(t_{k+1})\widehat{\phi}_{l_3}(t_k)\nn\\
&\;\;\;\;\; + (\widehat{\phi}_{l_2}(t_{k+1}))^{\ast}\widehat{\phi}_{l_1}(t_{k+1}) (\widehat{\phi}_{l_3}(t_k)-\widehat{\phi}_{l_3}(t_{k+1})),\label{eq:cksum}
\end{align}
where $c^\ast$ is the complex conjugate of $c$.
 For the step size in \eqref{eq:tauc1-nl}, we know from \eqref{eq:cfl} that
for  $C=\frac{2\alpha}{\sin(\alpha\pi)}$,
\begin{equation}\label{eq:Sbd}
\vert S_{n,l,l_1,l_2,l_3}\vert \leq \frac{1}{\vert\sin(\tau \delta_{l,l_1,l_2,l_3}/2)\vert}\leq\frac{C}{\tau\vert\delta_{l,l_1,l_2,l_3}\vert},\quad \forall n\ge0.
\end{equation}
Combining \eqref{eq:rest}, \eqref{eq:lambdasum}, \eqref{eq:cksum} and \eqref{eq:Sbd}, we have
\begin{align}
\left\vert\sum_{k=0}^n\Lambda_{k,l,l_1,l_2,l_3}\right\vert
\lesssim & \  \tau^2\vert\delta_{l,l_1,l_2,l_3}\vert\sum\limits_{k=0}^{n-1}\bigg(
\left\vert\widehat{\phi}_{l_1}(t_k)-\widehat{\phi}_{l_1}(t_{k+1})\right\vert\left\vert\widehat{\phi}_{l_2}(t_k)\right\vert \left\vert\widehat{\phi}_{l_3}(t_k)\right\vert\nonumber\\
&\ +\left\vert\widehat{\phi}_{l_1}(t_{k+1})\right\vert\left\vert\widehat{\phi}_{l_2}(t_k)-\widehat{\phi}_{l_2}(t_{k+1})\right\vert \left\vert\widehat{\phi}_{l_3}(t_k)\right\vert \nonumber\\
&\ +\left\vert\widehat{\phi}_{l_1}(t_{k+1})\right\vert\left\vert\widehat{\phi}_{l_2}(t_{k+1})\right\vert \left\vert\widehat{\phi}_{l_3}(t_k)-\widehat{\phi}_{l_3}(t_{k+1})\right\vert\bigg)\nonumber\\
&\ + \tau^2\vert\delta_{l,l_1,l_2,l_3}\vert \left\vert\widehat{\phi}_{l_1}(t_n)\right\vert\left\vert\widehat{\phi}_{l_2}(t_n)\right\vert \left\vert\widehat{\phi}_{l_3}(t_n)\right\vert.
\label{eq:sumlambda}
\end{align}
 Following the discussions in the proof of Theorem \ref{thm:im_linear}, for the non-resonance step size in \eqref{eq:tauc2-nl},
the same bound in \eqref{eq:sumlambda} holds by replacing $\vert\delta_{l,l_1,l_2,l_3}\vert$ with $\vert\delta_{l,l_1,l_2,l_3}\vert^{2+\nu_3}$. The rest arguments are almost the same as those for the \eqref{eq:tauc1-nl} case and we shall only treat the step size \eqref{eq:tauc1-nl} below.

For $l\in\mathcal{T}_{N_0}$ and $(l_1,l_2,l_3)\in\mathcal{I}_l^{N_0}$, there holds
\begin{equation}
(1+\vert\mu_l\vert)\vert\delta_{l,l_1,l_2,l_3}\vert \leq (1+\vert\mu_l\vert) \bigg[(\sum\limits_{j=1}^3\mu_{l_j})^2 +\sum_{j=1}^3\mu_{l_j}^2\bigg] \lesssim \prod_{j=1}^3(1+\mu^2_{l_j})^{3/2}.
\label{eq:mlbd}
\end{equation}
Based on \eqref{eq:remainder-dec}, \eqref{eq:sumlambda} and \eqref{eq:mlbd},  we have from \eqref{eq:final},
\begin{align}
&\|\mathcal{L}^n \|^2_{H^1}  \nn \\
& = \ \varepsilon^4
\sum\limits_{l\in\mathcal{T}_{N_0}}\left(1+\mu^2_l\right) \big\vert\sum\limits_{(l_1,l_2,l_3)\in\mathcal{I}_l^{N_0}}\sum\limits_{k=0}^n\Lambda_{k,l,l_1,l_2,l_3}\big\vert^2 \nn \\
& \lesssim \ \varepsilon^4\tau^4
\bigg\{\sum_{l\in\mathcal{T}_{N_0}}\bigg(\sum\limits_{(l_1,l_2,l_3)\in\mathcal{I}_l^{N_0}}\left\vert\widehat{\phi}_{l_1}(t_n)\right\vert\left\vert\widehat{\phi}_{l_2}(t_n)\right\vert \left\vert\widehat{\phi}_{l_3}(t_n)\right\vert\prod_{j=1}^3(1+\mu_{l_j}^2)^{\frac{3}{2}} \bigg)^2
\nn \\
&\quad+n \sum\limits_{k=0}^{n-1}
\sum_{l\in\mathcal{T}_{N_0}}\bigg[\bigg(\sum\limits_{(l_1,l_2,l_3)\in\mathcal{I}_l^{N_0}}
\left\vert\widehat{\phi}_{l_1}(t_k)-\widehat{\phi}_{l_1}(t_{k+1})\right\vert\left\vert\widehat{\phi}_{l_2}(t_k)\right\vert \left\vert\widehat{\phi}_{l_3}(t_k)\right\vert \prod_{j=1}^3(1+\mu_{l_j}^2)^{\frac{3}{2}} \bigg)^2\nn\\
& \;\;\;\; +\bigg(\sum\limits_{(l_1,l_2,l_3)\in\mathcal{I}_l^{N_0}}
\left\vert\widehat{\phi}_{l_1}(t_{k+1})\right\vert\left\vert\widehat{\phi}_{l_2}(t_k)-\widehat{\phi}_{l_2}(t_{k+1})\right\vert \left\vert \widehat{\phi}_{l_3}(t_k)\right\vert\prod_{j=1}^3(1+\mu_{l_j}^2)^{\frac{3}{2}} \bigg)^2\nn\\
& \;\;\;\; +\bigg(\sum\limits_{(l_1,l_2,l_3)\in\mathcal{I}_l^{N_0}}
\left\vert\widehat{\phi}_{l_1}(t_{k+1})\right\vert\left\vert\widehat{\phi}_{l_2}(t_{k+1})\right\vert \left\vert\widehat{\phi}_{l_3}(t_k)-\widehat{\phi}_{l_3}(t_{k+1})\right\vert\prod_{j=1}^3(1+\mu_{l_j}^2)^{\frac{3}{2}} \bigg)^2\bigg]\bigg\}.
\label{eq:sumlambda-2}
\end{align}
Introducing the auxiliary function $\xi(x)=\sum_{l\in\mathbb{Z}}(1+\mu_l^2)^{\frac{3}{2}}\left\vert\widehat{\phi}_l(t_n)\right\vert e^{i\mu_l(x-a)}$, where $\xi(x)\in H_{\rm per}^{m-3}(\Omega)$ implied by assumption (D) and $\|\xi\|_{H^s}\lesssim \|\psi(t_n)\|_{H^{s+3}}$. Expanding $\vert\xi(x)\vert^2\xi(x)=\sum\limits_{l\in\mathbb{Z}}\sum\limits_{l_1-l_2+l_3=l} \prod_{j=1}^3\left((1+\mu_{l_j}^2)^{\frac{3}{2}}\left\vert\widehat{\phi}_{l_j}(t_n)\right\vert\right) e^{i\mu_l(x-a)}$, we could obtain that
\begin{align}
&\sum_{l\in\mathcal{T}_{N_0}} \bigg(\sum\limits_{(l_1,l_2,l_3)\in\mathcal{I}_l^{N_0}}\left\vert\widehat{\phi}_{l_1}(t_n)\right\vert\left\vert\widehat{\phi}_{l_2}(t_n)\right\vert \left\vert\widehat{\phi}_{l_3}(t_n)\right\vert\prod_{j=1}^3(1+\mu_{l_j}^2)^{\frac{3}{2}} \bigg)^2\nonumber\\
& \leq \left\|\vert\xi(x)\vert^2\xi(x)\right\|^2_{L^2}\lesssim \left\|\xi(x)\right\|_{H^1}^6\lesssim \left\|\psi(t_k)\right\|_{H^4}^6 \lesssim1.
\end{align}
Noticing \eqref{eq:nl_twist}, we can estimate each terms in \eqref{eq:sumlambda-2} accordingly as
\begin{align}
&\left\| \sum\limits_{k=0}^n e^{-i(k+1)\tau\Delta}P_{N_0}\mathcal{J}_1(e^{i t_k\Delta}(P_{N_0}\phi(t_k)))\right\|^2_{H^1} \nonumber \\
& \lesssim \varepsilon^4\tau^4 \bigg[\left\|\phi(t_k)\right\|_{H^4}^6+n\sum\limits_{k=0}^{n-1}
\left\|\phi(t_k) - \phi(t_{k+1})\right\|_{H^4}^2(\left\|\phi(t_k)\right\|_{H^4} + \left\|\phi(t_{k+1})\right\|_{H^4})^4\bigg]\nn \\
&\lesssim  \varepsilon^4\tau^4+n^2\varepsilon^4\tau^4 (\varepsilon^2\tau)^2
\lesssim \varepsilon^4\tau^4,\quad n\leq q,
\label{eq:est-l2}
\end{align}
and \eqref{eq:final2} implies
\begin{align}
\|e^{n+1}\|_{H^1} \lesssim & \ h^{m-1} + \tau_0^{m-1} +\eps^2\tau^2+\eps^2\tau\sum_{k=0}^n\|e^k\|_{H^1},\quad 0\leq n\leq q. \label{eq:final3}
\end{align}
Using discrete Gronwall's inequality, we have
\begin{equation}
\|e^{q+1}\|_{H^1}\lesssim h^{m-1} + \eps^2\tau^2 + \tau_0^{m-1},\quad 0\leq q\leq\frac{T/\varepsilon^2}{\tau}-1,
\end{equation}
which implies the first inequality in \eqref{eq:error_nl} at $n = q+1$. There exists $h_0 >0$ and $\tau_0>0$, when $h \leq h_0$,  $\tau$ satisfies \eqref{eq:tauc1-nl} or \eqref{eq:tauc2-nl} with $0 < \tau \leq \tau_1/\varepsilon $  ($\tau_1$ sufficiently small), the triangle inequality yields that
\begin{equation*}
\left\|I_N\psi^{q+1}\right\|_{H^1} \leq \left\|\psi(x, t_{q+1})\right\|_{H^1} + \left\|e^{q+1}\right\|_{H^1} \leq M + 1, \quad 0 \leq q \leq \frac{T/\varepsilon^2}{\tau}-1,	
\end{equation*}
which means that the induction process for \eqref{eq:error_nl} is completed.  $\hfill\Box$

\begin{remark}
The improved uniform error bound for the NLSE in Theorem \ref{thm:eb_nl} is for the cubic nonlinearity without the external potential. It is straightforward to extend to the NLSE with the general nonlinearity $\eps^{2p}|u|^{2p}u$ ($p\in\mathbb{Z}^+$) and the external potential $\eps^{2p} V(x)$. The long-time dynamics of the NLSE with $O(\eps^{2p})$-nonlinearity and $O(1)$-initial data is equivalent to the NLSE with $O(1)$-nonlinearity and $O(\eps)$-initial data. The amplitude of the potential is also $O(\eps^{2p})$, where the scaling is to be consistent with the life-span of the NLSE. The improved $H^1$-error bound of the TSFP method for the NLSE with $\eps^{2p}|u|^{2p}u$ nonlinearity up to the time at $O(1/\eps^{2p})$ is $O(h^{m-1} + \eps^{2p}\tau^2 + \tau_0^{m-1})$. The discussions on removing the parameter $\tau_0$ in Remark \ref{rmk:tau0} could be adapted here and we omit the details for brevity.
\end{remark}

\subsection{Numerical results}
In this subsection, we present some numerical examples for the NLSE with $O(\varepsilon^2)$-nonlinearity in 1D and 2D to confirm the improved uniform error bound in $H^1$-norm.

First, we show the long-time temporal errors of the TSFP \eqref{eq:TSFP_nl} for the NLSE \eqref{eq:nl_1D} on 1D domain $[0,2\pi]$. The initial data is chosen as
\begin{equation}
\psi_0(x) = 2/(2+\sin^2(x)), \quad x \in [0, 2\pi].	
\end{equation}

\begin{figure}[ht!]
\centerline{\includegraphics[width=12cm,height=5.5cm]{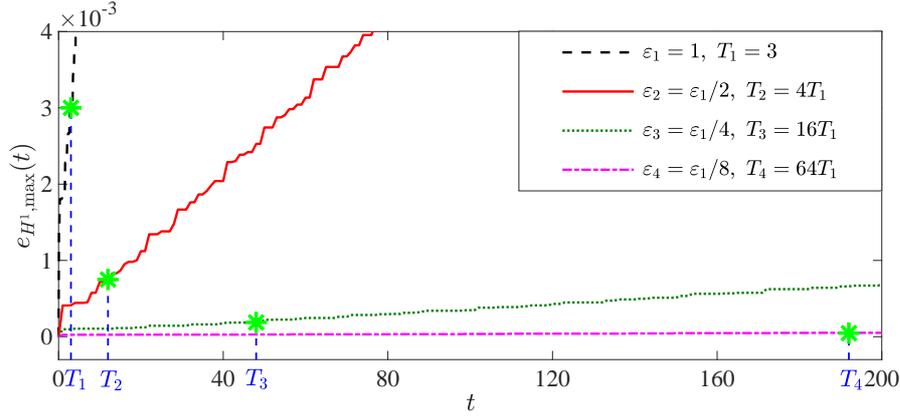}}
\caption{Long-time temporal errors in $H^1$- norm of the TSFP \eqref{eq:TSFP_nl} for the NLSE \eqref{eq:nl_1D} with different $\varepsilon$.}
\label{fig:nonlinear_long}
\end{figure}

\begin{figure}[ht!]
\begin{minipage}{0.49\textwidth}
\centerline{\includegraphics[width=6.3cm,height=5.5cm]{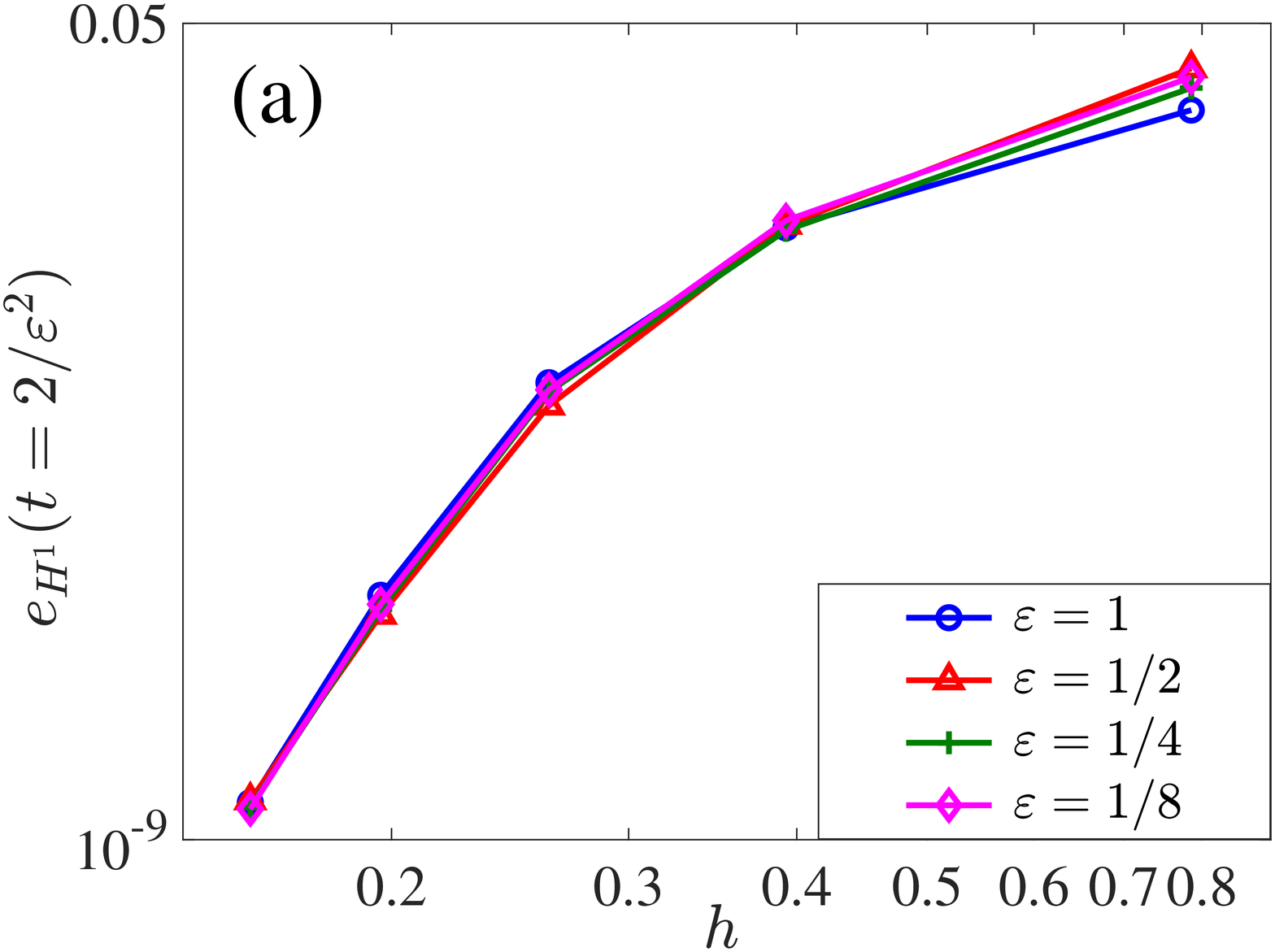}}
\end{minipage}
\begin{minipage}{0.49\textwidth}
\centerline{\includegraphics[width=6.3cm,height=5.5cm]{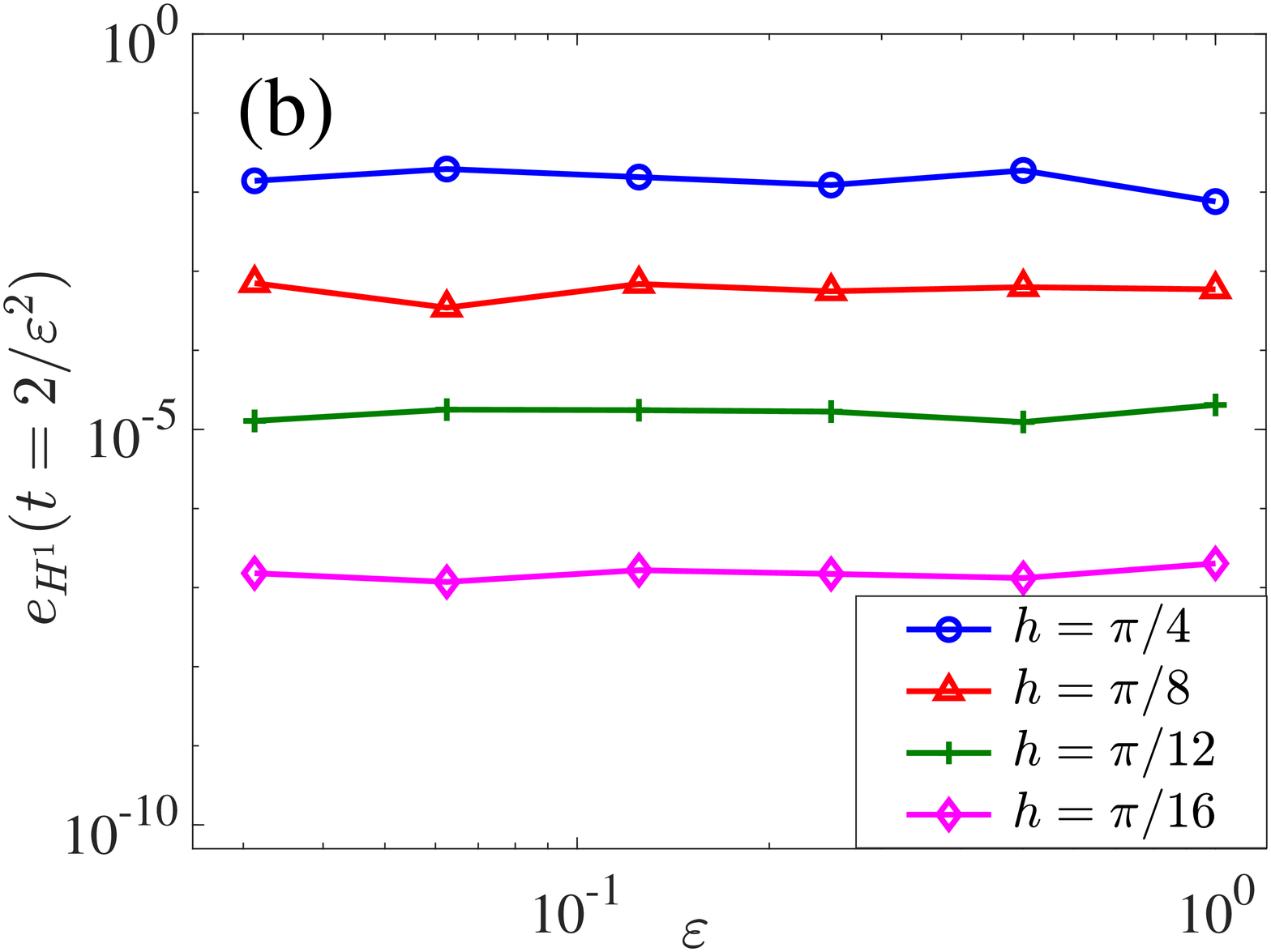}}
\end{minipage}
\caption{Long-time spatial errors in $H^1$-norm of the TSFP \eqref{eq:TSFP_nl} for the NLSE in \eqref{eq:nl_1D} at $t = 2/\varepsilon^2$.}
\label{fig:nonlinear_spatial}
\end{figure}

\begin{figure}[ht!]
\begin{minipage}{0.49\textwidth}
\centerline{\includegraphics[width=6.3cm,height=5.5cm]{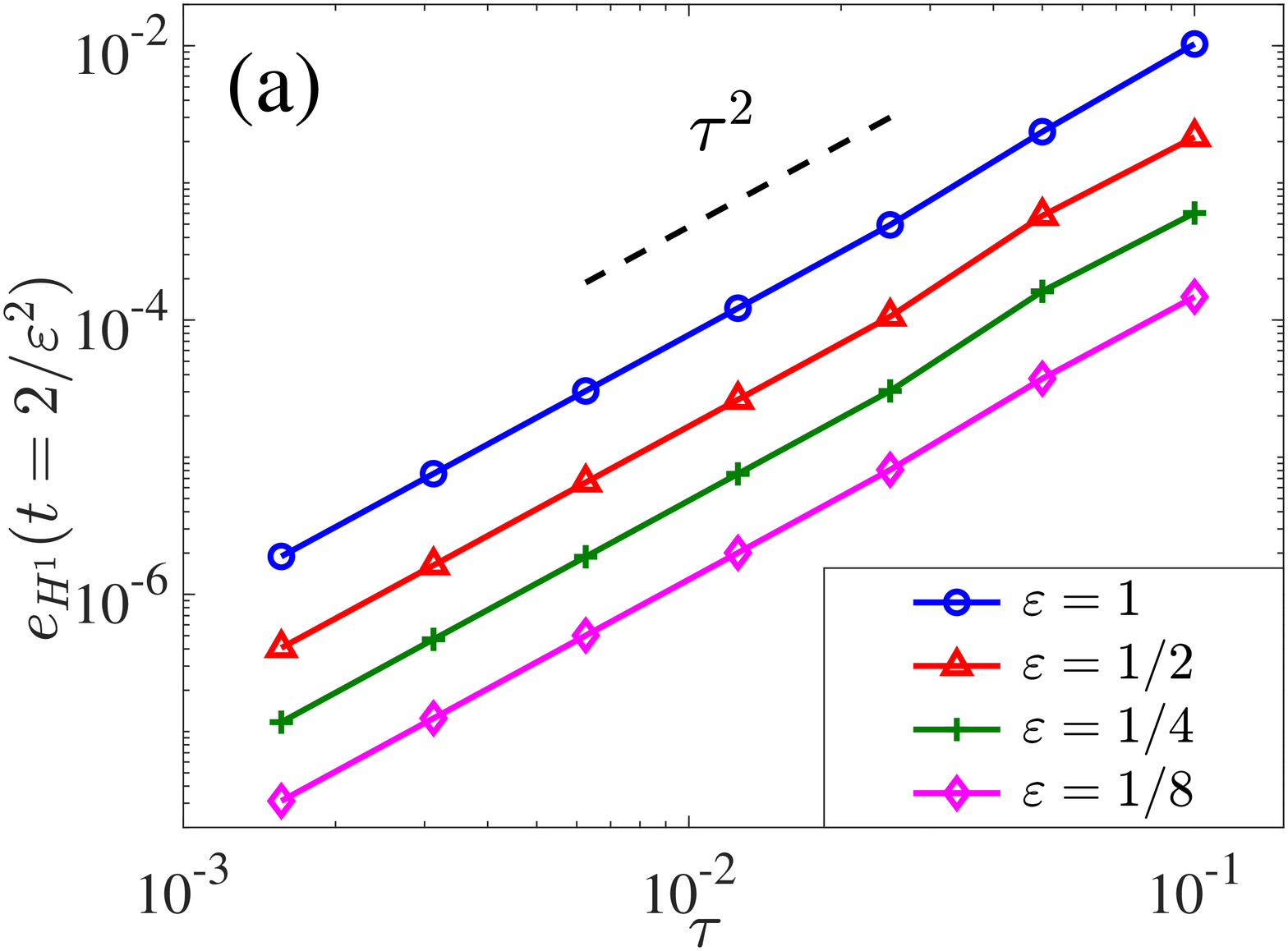}}
\end{minipage}
\begin{minipage}{0.49\textwidth}
\centerline{\includegraphics[width=6.3cm,height=5.5cm]{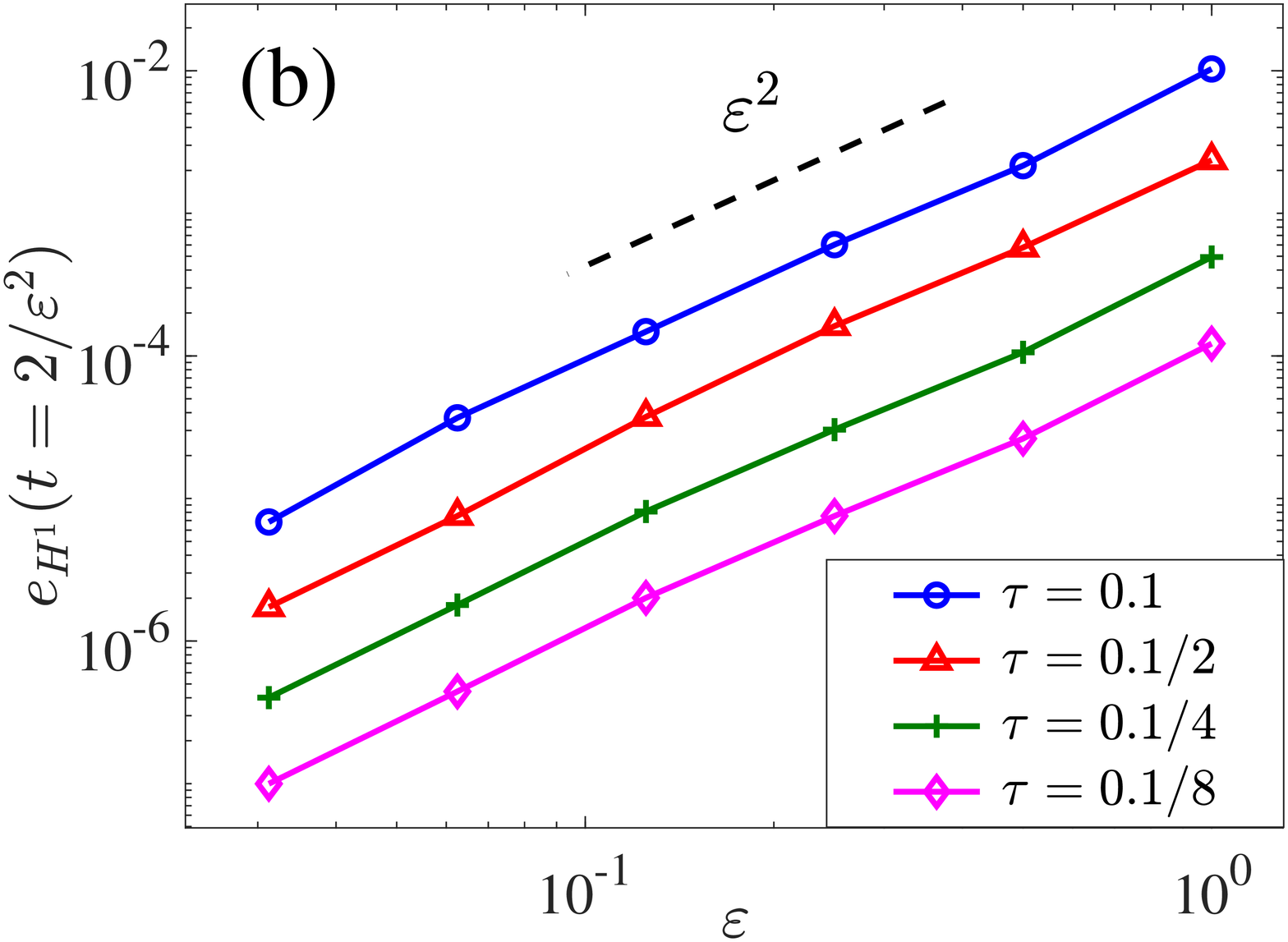}}
\end{minipage}
\caption{Long-time temporal errors in $H^1$-norm of the TSFP \eqref{eq:TSFP_nl} for the NLSE \eqref{eq:nl_1D} at $t = 2/\varepsilon^2$.}
\label{fig:nonlinear_temporal}
\end{figure}

Figure \ref{fig:nonlinear_long}  plots the long-time errors in $H^1$-norm of the TSFP method for the NLSE with a fixed time step $\tau$ and different $\varepsilon$, which indicates that the global errors in $H^1$-norm behave like $O(\eps^2\tau^2)$ up to the $O(1/\eps^2)$ time. Then, we show the spatial and temporal errors of the TSFP \eqref{eq:TSFP_nl} for the NLSE \eqref{eq:nl_1D}. Figure \ref{fig:nonlinear_spatial} \& Figure \ref{fig:nonlinear_temporal} depict the long-time spatial and temporal errors of the TSFP \eqref{eq:TSFP_nl} for the NLSE \eqref{eq:nl_1D} at $t= 2/\varepsilon^2$, respectively. Similar to the linear case, Figure \ref{fig:nonlinear_spatial} shows the spectral accuracy of the TSFP method for the NLSE in space and the spatial errors are independent of the small parameter $\varepsilon$. Each line in Figure \ref{fig:nonlinear_temporal} (a) corresponds to a fixed $\varepsilon$ and shows the global errors in $H^1$-norm versus the time step $\tau$, which confirms the second-order convergence of the TSFP method in time. Figure \ref{fig:nonlinear_temporal} (b) again validates that the global errors in $H^1$-norm behave like $O(\eps^2\tau^2)$ up to the $O(1/\eps^2)$ time.

\begin{figure}[ht!]
\begin{minipage}{0.49\textwidth}
\centerline{\includegraphics[width=6.3cm,height=5.5cm]{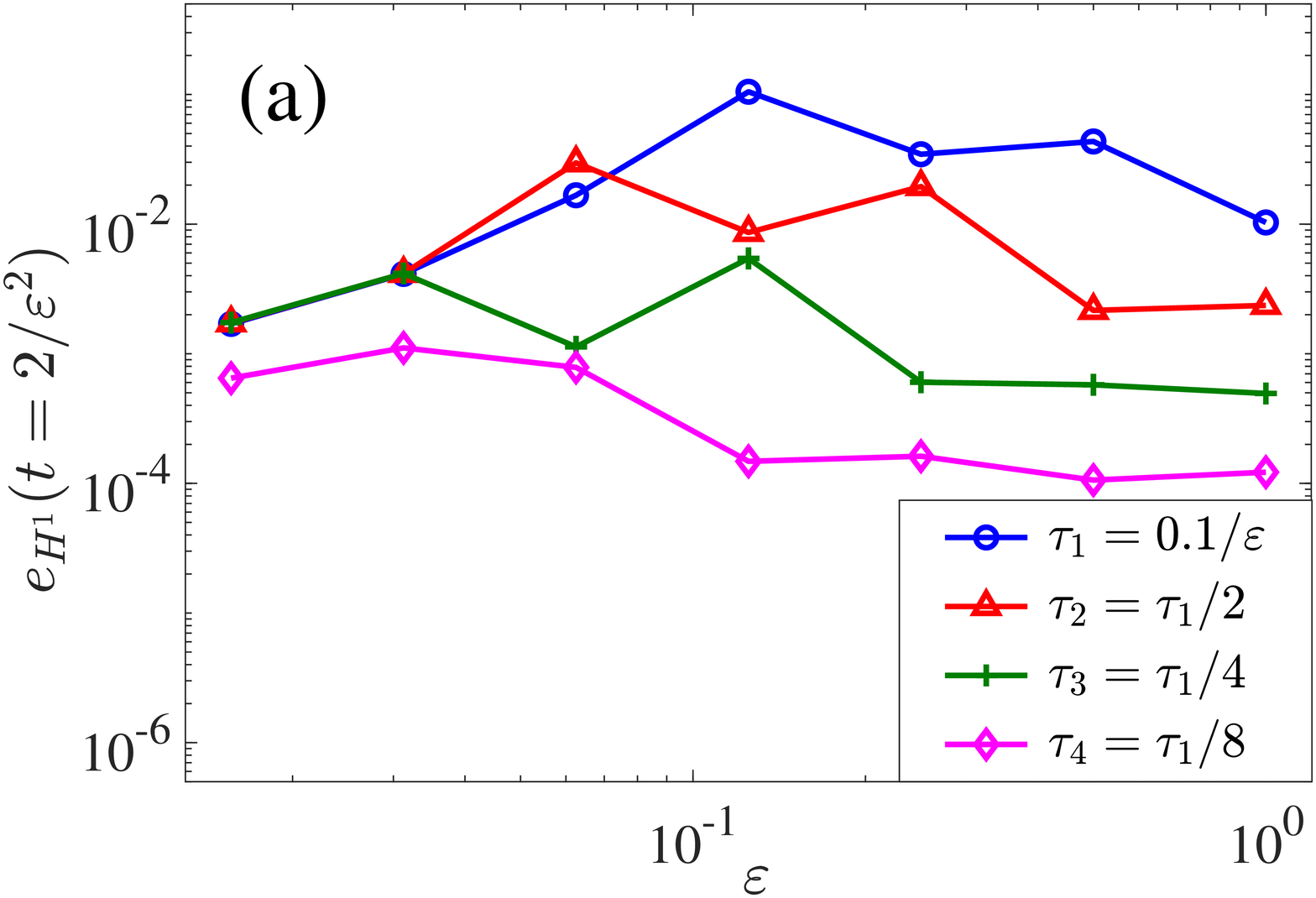}}
\end{minipage}
\begin{minipage}{0.49\textwidth}
\centerline{\includegraphics[width=6.3cm,height=5.5cm]{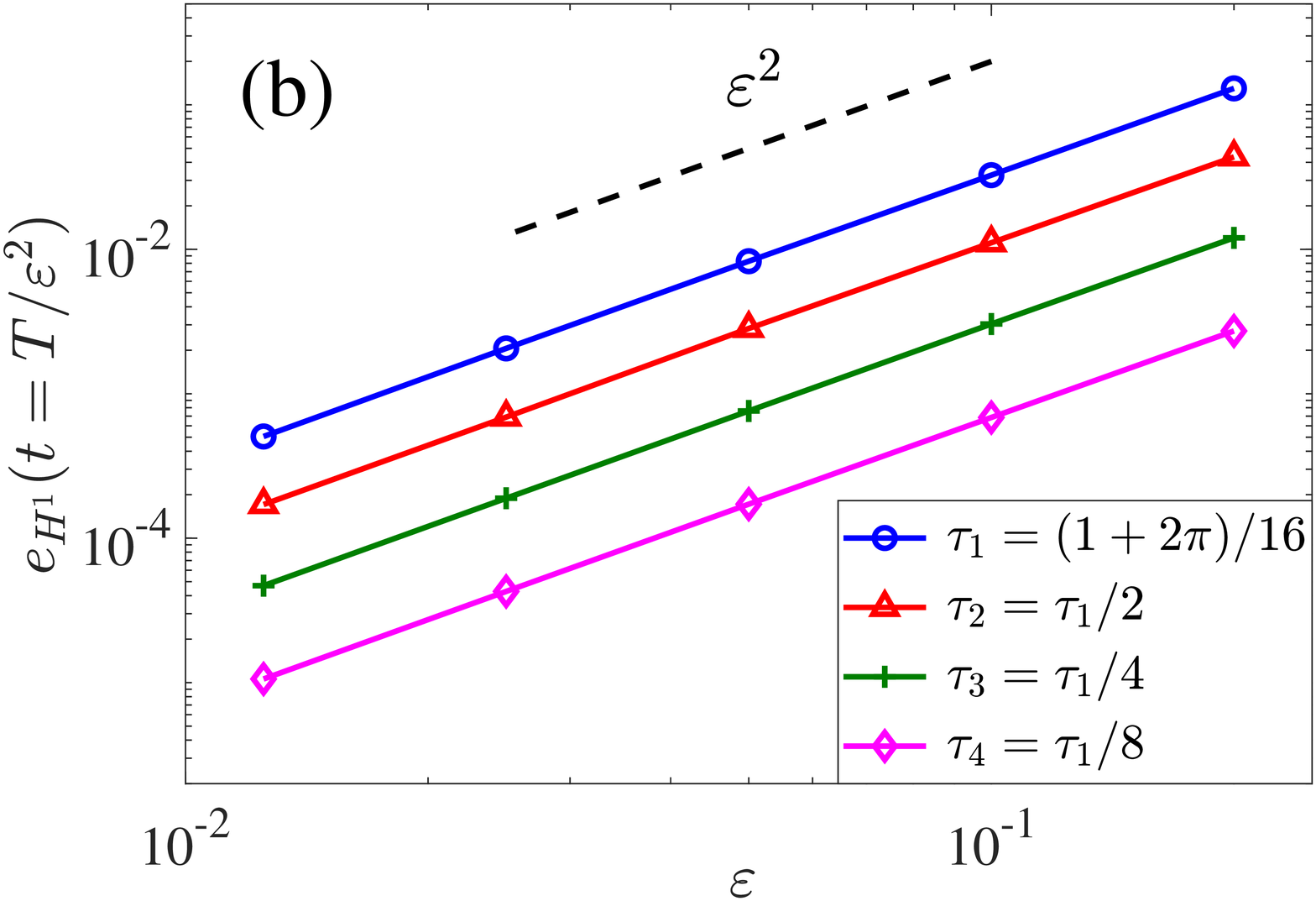}}
\end{minipage}
\caption{Long-time temporal errors in $H^1$-norm of the TSFP \eqref{eq:TSFP} for the NLSE \eqref{eq:nl_1D} with large time step size.}
\label{fig:nl_large}
\end{figure}

Figure \ref{fig:nl_large} displays the long-time errors of the TSFP method for the NLSE  \eqref{eq:nl_1D} with large time step size. Each line in Figure \ref{fig:nl_large} (a) plots the long-time errors with  $\tau = O(1/\eps)$, which are almost  constants for different $\eps$, confirming  the error bounds. In Figure \ref{fig:nl_large} (b), we choose the larger time step size  satisfying the non-resonance condition, which demonstrates the improved uniform error bounds with large time step size.

Then, we show an example in 2D with the irrational aspect ratio of the domain. We choose the domain $\Omega=(0, \pi) \times (0, 1)$  and the initial data
\begin{equation}
\psi_0(x, y) = \frac{1}{1+ \sin^2(2x)} + \sin(2\pi y), \quad {\bf x}=(x , y) \in [0, \pi] \times [0, 1].
\label{eq:initial_2D}	
\end{equation}

\begin{figure}[ht!]
\centerline{\includegraphics[width=12cm,height=5.5cm]{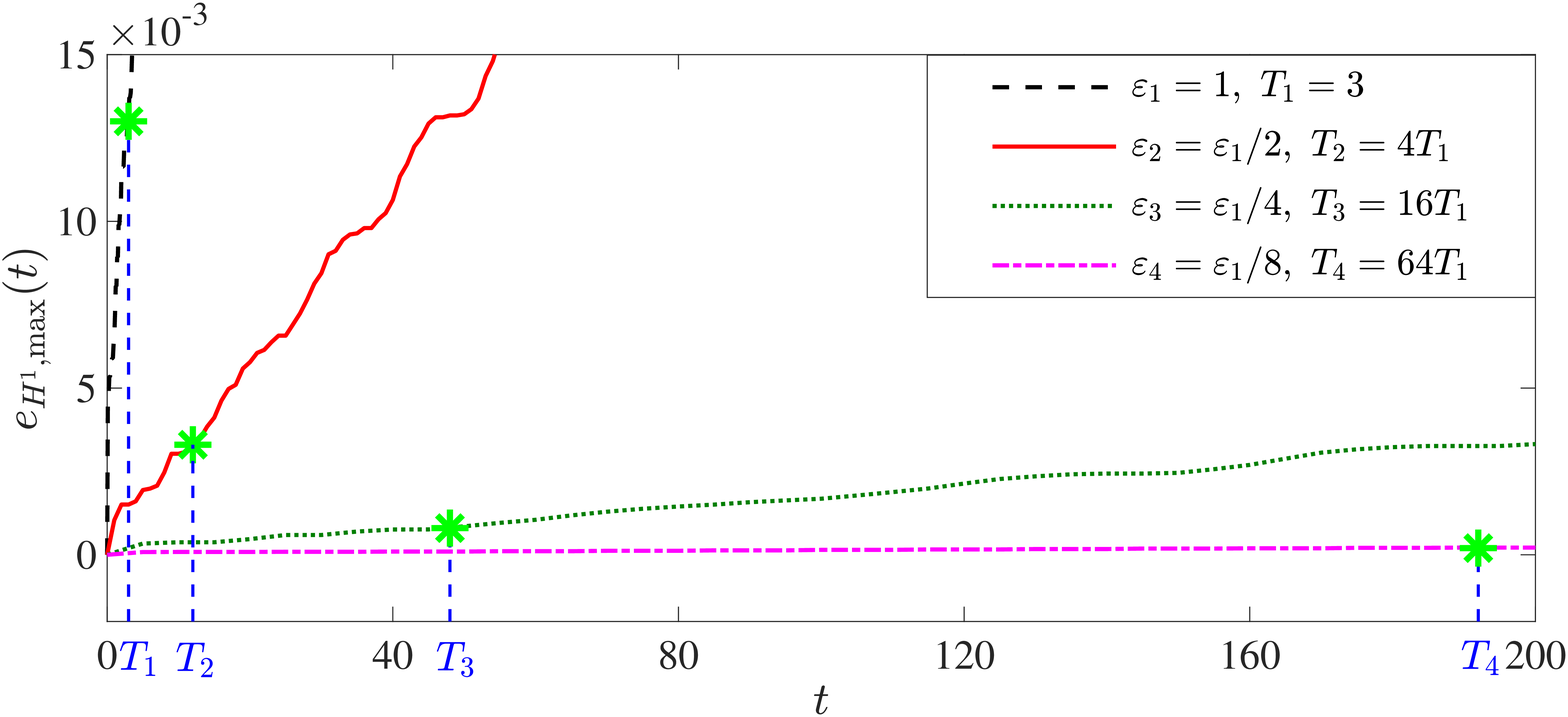}}
\caption{Long-time temporal errors in $H^1$-norm of the TSFP method for the NLSE \eqref{eq:WNE} in 2D with different $\eps$.}
\label{fig:2D_nonlinear_long}
\end{figure}

\begin{figure}[ht!]
\begin{minipage}{0.49\textwidth}
\centerline{\includegraphics[width=6.3cm,height=5.5cm]{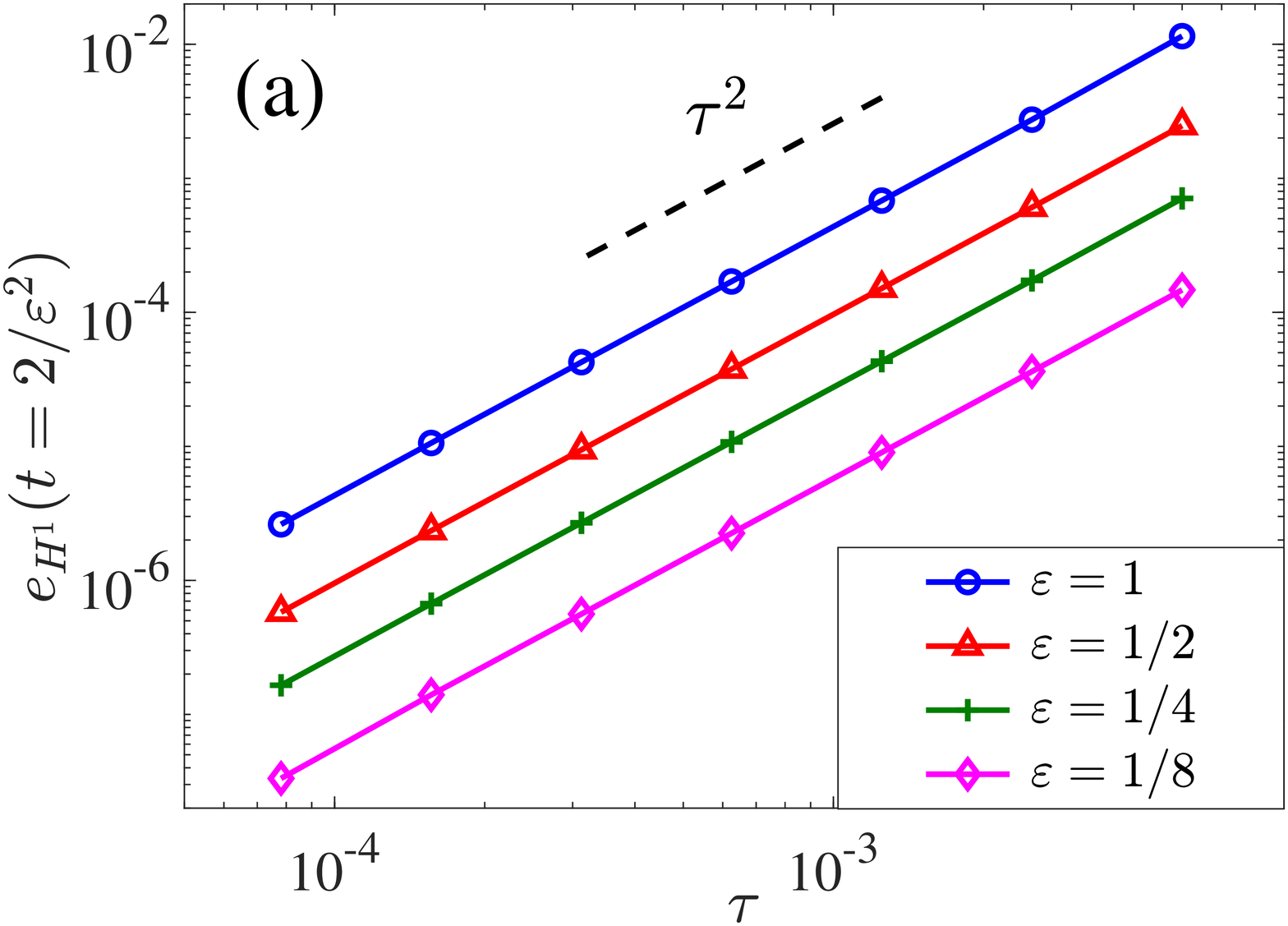}}
\end{minipage}
\begin{minipage}{0.49\textwidth}
\centerline{\includegraphics[width=6.3cm,height=5.5cm]{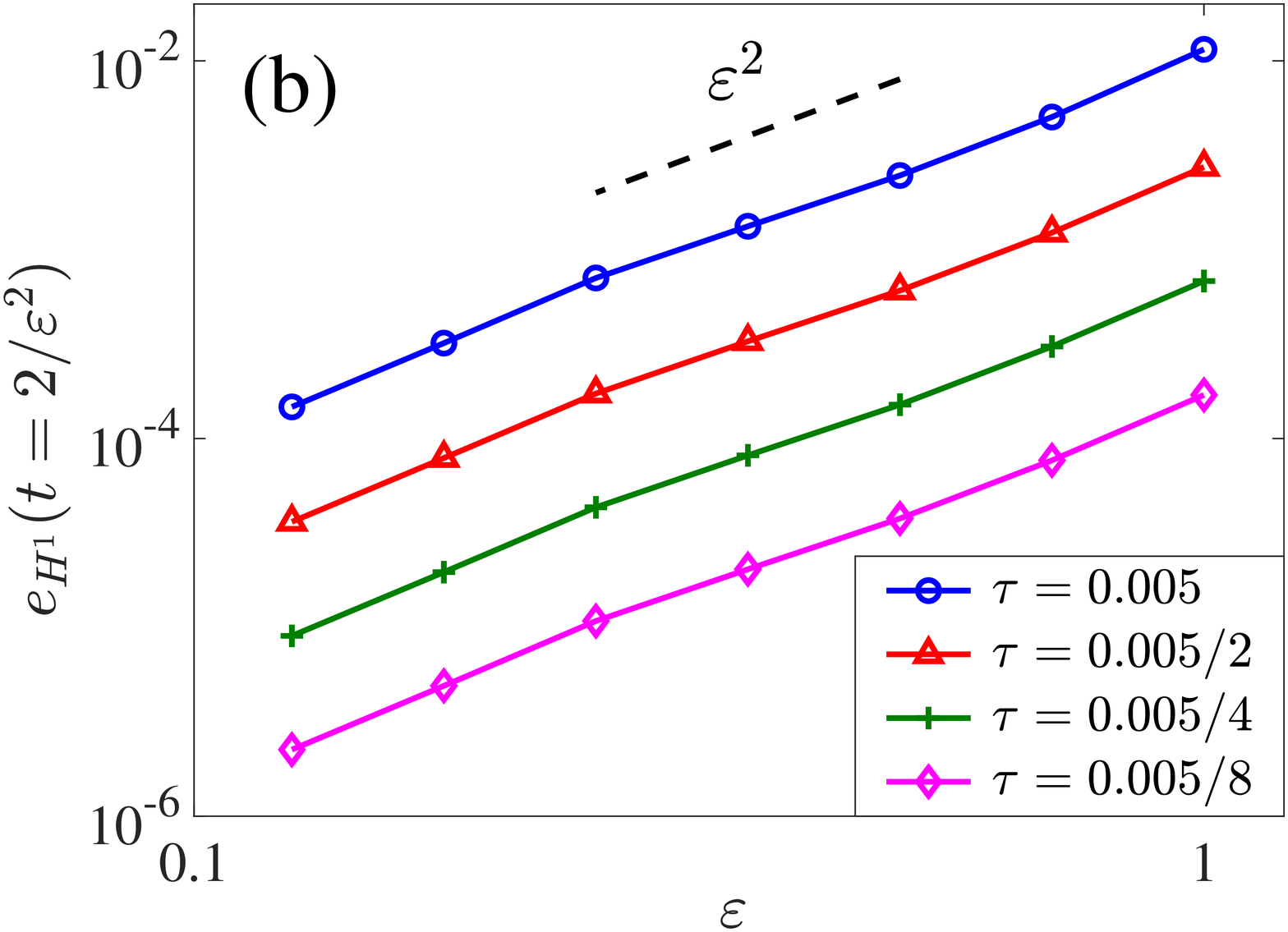}}
\end{minipage}
\caption{Long-time temporal errors in $H^1$-norm of the TSFP method for the NLSE \eqref{eq:WNE} in 2D at $t = 1/\varepsilon^2$.}
\label{fig:2D_temporal}
\end{figure}

 Figure \ref{fig:2D_nonlinear_long} plots the long-time temporal errors in $H^1$-norm of the TSFP method for the NLSE in 2D with a fixed time step $\tau$ and different $\varepsilon$, which confirms that the improved uniform error bound in $H^1$-norm at $O(\eps^2\tau^2)$ up to the  $O(1/\eps^2)$ time is also suitable for the irrational aspect ratio of the domain. Figure \ref{fig:2D_temporal} depicts the long-time errors for the TSFP method for the NLSE in 2D  at $t = 1/\eps^2$, which again indicates that the TSFP method is second-order in time and validates the improved uniform error bound in $H^1$-norm up to the time at $O(1/\eps^2)$.


\section{Conclusions}
Improved uniform error bounds for the time-splitting Fourier pseudospectral (TSFP) methods for the long-time dynamics of the Schr\"odinger equation with small potential and the nonlinear Schr\"odinger equation (NLSE) with weak nonlinearity were rigorously established.
 For the Schr\"odinger equation with small potential, the linear growth of the uniform error bound in $L^2$-norm for the TSFP method was strictly proven with the aid of the unitary property of the solution flow in $L^2(\Omega)$. By introducing a new technique of regularity compensation oscillation (RCO), the improved uniform error bound in $H^1$-norm was carried out at $O(h^{m-1}+\eps\tau^2)$ up to the $O(1/\eps)$ time.  In addition, the RCO technique was extended to show the improved uniform error bound $O(h^{m-1}+\eps^2\tau^2)$ for the TSFP method applied to the cubic NLSE with $O(\eps^2)$-nonlinearity up to the $O(1/\eps^2)$ time. Numerical results were presented to validate our error estimates and demonstrate that they are sharp. We remark here that the RCO technique has been adapted to establish improved uniform error bounds on time-splitting methods for the long-time dynamics of dispersive PDEs including the nonlinear Klein-Gordon equation \cite{BCF} and the nonlinear Dirac equation \cite{BCF1}.

 \section*{Acknowledgement}
The authors would like to thank the anonymous referee for the invaluable comments and suggestions.
\bibliographystyle{amsplain}

\begin{thebibliography}{}
%
%

         \bibitem{AK}
         G. D. Akrivis, \textit{Finite difference discretization of the cubic Schr\"odinger equation}, IMA J. Numer. Anal. \textbf{13} (1993), no. 1, 115--124.

        \bibitem{ABB}
        X. Antoine, W. Bao, and C. Besse, \textit{Computational methods for the dynamics of the nonlinear Schr\"odinger/Gross--Pitaevskii equations}, Comput. Phys. Commun. \textbf{184} (2013), no. 12, 2621--2633.

        \bibitem{BC}
        W. Bao and Y. Cai, \textit{Mathematical theory and numerical methods for Bose-Einstein condensation}, Kinet. Relat. Models \textbf{6} (2013), no. 1, 1--135.

           \bibitem{BC14}
        W. Bao and Y. Cai, \textit{Uniform and optimal error estimates of an exponential wave integrator sine pseudospectral method for the nonlinear Schr\"odinger equation with wave operator},  SIAM J. Numer. Anal. \textbf{52} (2014), no. 3, 1103--1127.

\bibitem{BCF}
W. Bao, Y. Cai and Y. Feng,
\textit{Improved uniform error bounds on time-splitting methods for long-time dynamics of the nonlinear Klein-Gordon equation with weak nonlinearity}, SIAM J. Numer. Anal., to appear.

\bibitem{BCF1}
W. Bao, Y. Cai and Y. Feng,
\textit{Improved uniform error bounds on time-splitting methods for the long-time dynamics of the weakly nonlinear Dirac equation}, arXiv: 2203.05886.

\bibitem{BCY20}
W. Bao, Y. Cai and J. Yin, \textit{Super-resolution of time-splitting methods for the Dirac equation in the nonrelativistic regime}, Math. Comp.  \textbf{89} (2020), 2141-2173.

\bibitem{BCY21}
W. Bao, Y. Cai and J. Yin, \textit{Uniform error bounds of time-splitting methods for the nonlinear Dirac equation in the nonrelativistic regime without magnetic potential}, SIAM J. Numer. Anal. \textbf{59} (2021), 1040-1066.

        \bibitem{BJP}
        W. Bao, D. Jaksch, and P. A. Markowish, \textit{Numerical solution of the Gross--Pitaevskii equation for Bose--Einstein condensation}, J. Comput. Phys. \textbf{187} (2003), no. 1, 318--342.

        \bibitem{BJP2}
        W. Bao, S. Jin, and P. A. Markowish, \textit{On time-splitting spectral approximations for the Schr\"odinger equation in the semiclassical regime}, J. Comput. Phys. \textbf{175} (2002), no. 2, 487--524.


         \bibitem{BS}
        W. Bao and J. Shen, \textit{A fourth-order time-splitting Laguerre--Hermite pseudospectral method for Bose--Einstein condensates}, SIAM J. Sci. Comput. \textbf{26} (2005), no. 6, 2010--2028.

        \bibitem{BBD}
        C. Besse, B. Bid\'egaray, and S. Descombes, \textit{Order estimates in time of splitting methods for the nonlinear Schr\"odinger equation}, SIAM J. Numer. Anal. \textbf{40} (2002), no. 1, 26--40.

        \bibitem{Bour}
        J. Bourgain, \textit{Fourier transform restriction phenomena for certain lattice subsets and application to nonlinear evolution equations. Part I: Schr\"odinger equations}, Geom. Funct. Anal. \textbf{3} (1993), no. 2, 107--156.

      \bibitem{Bour2}
      J. Bourgain, \textit{Growth of Sobolev norms in linear Schr\"odinger equations with quasi-periodic potential}, Comm. Math. Phys. \textbf{204} (1999), no. 1, 207--247.

        \bibitem{BGHS}
        T. Buckmaster, P. Germain, Z. Hani, and J. Shatah, \textit{Effective dynamics of the nonlinear Schr\"odinger equation on large domains}, Comm. Pure Appl. Math. \textbf{71} (2018), no. 7, 1407--1460.

       \bibitem{BGT}
       N. Burq, P. G\'erard, and N. Tzvetkov, \textit{Strichartz inequalities and the nonlinear Schr\"odinger equation on compact manifolds}, Amer. J. Math. \textbf{126} (2004), no. 3, 569--605.



         \bibitem{Carles}
         R. Carles, \textit{On Fourier time-splitting methods for nonlinear Schr\"odinger equations in the semiclassical limit}, SIAM J. Numer. Anal. \textbf{51} (2013), no. 6, 3232--3258.

         \bibitem{CCMM}
         F. Castella, P. Chartier, F. M\'ehats, and A. Murua, \textit{Stroboscopic averaging for the nonlinear Schr\"odinger equation}, Found. Comput. Math. \textbf{15} (2015), no. 2, 519--559.

        \bibitem{Caz}
        T. Cazenave, \textit{Semilinear Schr\"odinger Equations}, Courant Lect. Notes Math., 10, Amer. Math. Soc., Providence, RI, 2003.

         \bibitem{CCO}
        E. Celledoni, D. Cohen, and B. Owren, \textit{Symmetric exponential integrators with an application to the cubic Schr\"odinger equation}, Found. Comp. Math. \textbf{8} (2008), no. 3, 303--317.

	    \bibitem{CMTZ}
        P. Chartier, F. M{\'e}hats, M. Thalhammer, and Y. Zhang, \textit{Improved error estimates for splitting methods applied to highly-oscillatory nonlinear Schr{\"o}dinger equations}, Math. Comp. \textbf{85} (2016), no. 302, 2863--2885.

         \bibitem{CLL}
         W. Chen, X. Li, and D. Liang, \textit{Energy-conserved splitting FDTD Methods for Maxwell's equations}, Numer. Math. \textbf{108} (2008), no. 3, 445--485.

         \bibitem{CLL2}
         W. Chen, X. Li, and D. Liang, \textit{Energy-conserved splitting finite difference time domain methods for Maxwell's equations in three dimensions}, SIAM J. Numer. Anal.  \textbf{48} (2010), no. 4, 1530--1554.

         \bibitem{CHL}
         D. Cohen, E. Hairer, and C. Lubich, \textit{Modulated Fourier expansions of highly oscillatory differential equations}, Found. Comput. Math. \textbf{3} (2003), no. 4, 327--345.

         \bibitem{CK}
        J. Colliander, M. Keel, G. Staffilani, H. Takaoka, and T. Tao, \textit{Almost conservation laws and global rough solutions to a nonlinear Schr\"odinger equation}, Math. Res. Lett. \textbf{9} (2002), no. 5, 1--24.

\bibitem{Deb}
         A. Debussche and E. Faou, \textit{Modified energy for split-step methods applied to the linear Schr\"odinger equation}, SIAM J. Numer. Anal. \textbf{47} (2009), no. 5, 3705--3719.

         \bibitem{DFP}
         M. Delfour, M. Fortin, and G. Payre, \textit{Finite-difference solutions of a nonlinear Schr\"odinger equation}, J. Comput. Phys. \textbf{44} (1981), no. 2, 277--288.

    \bibitem{DG}
         G. Dujardin, \textit{Exponential Runge--Kutta methods for the Schr\"odinger equation}, Appl. Numer. Math. \textbf{59} (2009), no. 8, 1839--1857.

         \bibitem{DF}
         G. Dujardin and E. Faou, \textit{Normal form and long time analysis of splitting schemes for the linear Schr\"odinger equation with small potential}, Numer. Math. \textbf{108} (2007), no. 2, 223--262.


        \bibitem{ESY}
        L. Erd\text{\H o}s, B. Schlein, and H.-T. Yau, \textit{Derivation of the cubic non-linear Schr\"odinger equation from quantum dynamics of many-body systems}, Invent. Math. \textbf{167} (2007), no. 3, 515--614.


        \bibitem{Faou}
        E. Faou, \textit{Geometric Numerical Integration and Schr\"odinger Equations}, European Mathematical Society, Z\"urich, 2012.

         \bibitem{FGH}
         E. Faou, L. Gauckler, and Z. Hani, \textit{The weakly nonlinear large-box limit of the 2D cubic nonlinear Schr\"odinger equation}, J. Amer. Math. Soc. \textbf{29} (2016), no. 4, 915--982.

         \bibitem{FGL}
         E. Faou, L. Gauckler, and C. Lubich, \textit{Sobolev stability of plane wave solutions to the cubic nonlinear Schr\"odinger equation on a torus}, Comm. Partial Differential Equations \textbf{38} (2013), no. 7, 1123--1140.

         \bibitem{FGP}
         E. Faou, B. Gr\'ebert, and E. Paturel, \textit{Birkhoff normal form for splitting methods applied to semilinear Hamiltonian PDEs. I. Finite-dimensional discretization}, Numer. Math. \textbf{114} (2010), no. 3, 429--458.

         \bibitem{GL1}
         L. Gauckler and C. Lubich, \textit{Nonlinear Schr\"odinger equations and their spectral semi-discretizations over long times}, Found. Comput. Math. \textbf{10} (2010), no. 2, 141--169.

         \bibitem{GL2}
         L. Gauckler and C. Lubich, \textit{Splitting integrators for nonlinear Schr\"odinger equations over long times}, Found. Comput. Math. \textbf{10} (2010), no. 3, 275--302.

        \bibitem{HLW}
        E. Hairer, C. Lubich, and G. Wanner, \textit{Geometric Numerical Integration: Structure-Preserving Algorithms for Ordinary Differential Equations}, Springer, Berlin, 2002.


         \bibitem{HTT}
         S. Herr, D. Tataru, and N. Tzvetkov, \textit{Strichartz estimates for partially periodic solutions to Schr\"odinger equations in 4d and applications}, J. Reine Angew. Math.  \textbf{690} (2014), 65--78.

 	    \bibitem{HO}
	   M. Hochbruck and A. Ostermann, \textit{Exponential integrators}, Acta Numer. \textbf{19} (2010), 209--286.	

 \bibitem{JL00}
T. Jahnke and C. Lubich, \textit{Error bounds for exponential operator splittings}, BIT Numer. Math.  \textbf{40}(2000), 735--744.
       \bibitem{KAD}
       O. Karakashian, G. D. Akrivis, and V. A. Dougalis, \textit{On optimal order error estimates for the nonlinear Schr\"odinger equation}, SIAM J. Numer. Anal. \textbf{30} (1993), no. 2, 377--400.



       \bibitem{LU}
       C. Lubich, \textit{On splitting methods for Schr\"odinger-Poisson and cubic nonlinear Schr\"odinger equations}, Math. Comp. \textbf{77} (2008), no. 264, 2141--2153.


       \bibitem{MQ}
       R. I. McLachlan and G. R. W. Quispel, \textit{Splitting methods}, Acta Numer. \textbf{11} (2002), 341--434.


\bibitem{SZ}
   {Z. Shang},    \textit{Resonant and Diophantine step sizes in computing invariant tori of Hamiltonian systems}, Nonlinearity \textbf{13} (2000), 299--308.

       \bibitem{ST}
       {J. Shen,  T. Tang, and L. Wang}, \textit{Spectral Methods: Algorithms, Analysis and Applications}, Springer-Verlag Berlin Heidelberg, 2011.


         \bibitem{Strang}
         G. Strang, \textit{On the construction and comparison of difference schemes}, SIAM J. Numer. Anal. \textbf{5} (1968), no. 3, 506--517.

        \bibitem{SS}
        C. Sulem and P. Sulem, \textit{The Nonlinear Schr\"odinger Equation: Self-Focusing and Wave Collapse}, Springer, New York, 1999.


       	\bibitem{Tao}
        T. Tao, \textit{Nonlinear Dispersive Equations: Local and Global Analysis}, Amer. Math. Soc., Providence, RI, 2006.

        \bibitem{Thal}
        M. Thalhammer, \textit{High-order exponential operator splitting methods for time-dependent Schr\"odinger equations}, SIAM J. Numer. Anal. \textbf{46} (2008), no. 4, 2022--2038.

        \bibitem{Thal2}
        M. Thalhammer, \textit{Convergence analysis of high-order time-splitting pseudospectral methods for nonlinear Schr\"odinger equations}, SIAM J. Numer. Anal. \textbf{50} (2012), no. 6, 3231--3258.


         \bibitem{Wang}
         W.-M. Wang, \textit{Bounded Sobolev norms for linear Schr\"odinger equations under resonant perturbations}, J. Func. Anal. \textbf{254} (2008), no. 11, 2926--2946.

         \bibitem{WH}
         J. A. C. Weideman and B. M. Herbst, \textit{Split-step methods for the solution of the nonlinear Schr\"odinger equation}, SIAM J. Numer. Anal. \textbf{23} (1986), no. 3,  485--507.


\end{thebibliography}

\end{document}